\documentclass[12pt]{article}

\usepackage[english]{babel}

\usepackage[letterpaper,top=2cm,bottom=2cm,left=3cm,right=3cm,marginparwidth=1.75cm]{geometry}    

\usepackage{amsmath}
\usepackage{graphicx}
\usepackage[colorlinks=true, allcolors=blue]{hyperref}

\usepackage{subfigure}
\usepackage{amssymb,amsfonts,amsmath,latexsym,amsthm}

\usepackage[nottoc]{tocbibind}

\usepackage[square,numbers]{natbib}
\newcommand{\authorblock}[3]{%
  \begin{minipage}[t]{0.28\textwidth}
    \centering
    #1\\[0.35em]
    {\small #2}\\[0.4em]
    {\scriptsize\ttfamily #3}
  \end{minipage}%
}
\date{}

\newcommand{\R}{\mathbb{R}}

\newcommand{\E}{\mathbb{E}}
\newcommand{\ord}{m_0}

\newtheorem{theorem}{Theorem}
\newtheorem{corollary}{Corollary}
\newtheorem{proposition}{Proposition}
\newtheorem{lemma}{Lemma}

\newtheorem{definition}{Definition}

\newtheorem{remark}{Remark}
\newtheorem{example}{Example}

\providecommand{\keywords}[1]
{
  \small	
  \textbf{\textit{Keywords---}} #1
}

\title{Adaptive estimation in regression models for weakly dependent data and explanatory variables with unknown density}

\author{%
  \authorblock{Karine Bertin\footnotemark[1]}%
    {CIMFAV-INGEMAT\\Universidad de Valpara\'iso\\Valpara\'iso, Chile}%
  {karine.bertin@uv.cl}%
  \and
  \authorblock{Lisandro Ferm\'in\footnotemark[2]}%
    {Aix-Marseille University\\CNRS, AMSE\\Marseille, France}%
{lisandro.fermin@univ-amu.fr}%
  \and
  \authorblock{Miguel Padrino\footnotemark[3]}%
    {Departamento de Matem\'aticas\\Instituto Venezolano de Investigaciones Cient\'ificas\\Caracas, Venezuela}%
{padrinomiguel@gmail.com}%
}

\begin{document}
\maketitle

\begingroup
\renewcommand{\thefootnote}{\fnsymbol{footnote}}
\footnotetext[1]{Supported by FONDECYT regular grants 1221373
and 1230807 from ANID-Chile, MathAmSud AMSUD250012,
and the Centro de Modelamiento Matem\'atico (CMM) BASAL fund
FB210005 for centers of excellence from ANID-Chile.}
\footnotetext[2]{The project leading to this publication has received
funding from the French government under the ``France 2030''
investment plan managed by the French National Research Agency
(reference: ANR-17-EURE-0020) and from the Excellence Initiative
of Aix-Marseille University -- A*MIDEX.}
\footnotetext[3]{Supported by FONDECYT regular grant
1221373 from ANID-Chile.}
\endgroup

\begin{abstract}
This article is dedicated to the estimation of the regression function when the explanatory variable is a weakly dependent process whose correlation coefficient decays exponentially and has an unknown density function \( g \). The accuracy of the estimation is measured using pointwise risks. The regression function is expressed as the ratio between a function \( m \) and the density \( g \). A data-driven procedure based on kernel estimation is proposed to estimate the function \( m \). The bandwidth is selected using the approach of Goldenshluger and Lepski (GL); we prove that the resulting estimator of the function \( m \) satisfies an oracle-type inequality, which allows us to establish its adaptivity. According to \citet{bertin2017pointwise}, a GL estimator of the density function \( g \) also satisfies the aforementioned properties. Subsequently, a Nadaraya-Watson type estimator is used, similar to that proposed in \citet{fermin2017robbins}, with the difference that this estimator depends on two different bandwidths—one for the numerator and another for the denominator. Finally, it is shown that the proposed estimator for the regression function is adaptive with respect to the smoothness of $m$ and $g$. Additionally, unsupervised statistical learning techniques are described and applied to calibrate the method, and some simulations are provided to illustrate the performance of the proposed approach.
\end{abstract}


\tableofcontents

\section{Introduction}

In this article, we address the problem of adaptive nonparametric estimation using kernel estimators for regression functions in a univariate context and pointwise risk, based on observations
 \(\{(X_i, Y_i)\}_{i=1}^n\) from the regression model \(Y_i = r(X_i) + \epsilon_i\), where the explanatory variables \(\{X_i\}_{i=1}^n\) are weakly dependent and share an unknown common density \(g(\cdot)\). Adaptive estimation methods allow one to obtain estimators that converge at an optimal rate over a general class of functions, without requiring knowledge of the smoothness of the target function. Among the various adaptive estimation methods, the Goldenshluger–Lepski (GL) method stands out, introduced in \citet{goldenshluger2011bandwidth} and \citet{goldenshluger2014adaptive}.

Over the past decade, the GL method has been used in several statistical models—such as density estimation, regression, conditional density estimation, white noise models, among others—to derive oracle-type inequalities and, consequently, adaptive estimators. This is especially the case in independent settings (see, for example, \citet{BLR}, \citet{chichignoud2017adaptive}, \citet{goldenshluger2011bandwidth}, \citet{GL5}, \citet{Kouakou2020Lp}, \citet{Kouakou2023Sup}, \citet{lepski2014adaptive}). Furthermore, it has also been applied in certain dependent contexts, such as in \citet{asin2017adaptive}, \citet{bertin2017pointwise}, \citet{bertin2020adaptivebis}, \citet{bertin2020adaptive}, \citet{mpadrino2024adaptiregre}, and \citet{comte2017adaptive}.

Regarding adaptive results in regression models and dependence contexts, we also mention \citet{baraud2001model} for the \(L^2\)risk, who consider the case of 
$\beta$-mixing data and autoregressive models. They study least-squares estimators over families of models and show that a penalized least-squares estimator is adaptive. On the other hand, \citet{asin2017adaptive} consider adaptive nonparametric estimation of density and regression in independent and weakly dependent contexts, where the target function belongs to the Hilbert space \(L^2 := L^2[0, 1]\). The estimator is based on an orthogonal series approach, where the selection of the dimension parameter is fully data-driven, inspired by the GL method.

In \citet{mpadrino2024adaptiregre}, adaptive nonparametric kernel estimation of the bounded regression function \(r(\cdot)\) is carried out using the GL method under pointwise risk. This estimation is performed in a Hölder space of regularity \(\beta > 0\), for a regression model with a weakly dependent explanatory variable that is identically distributed with a known and bounded density \(g(\cdot)\). A Nadaraya–Watson–type estimator is used:
\begin{equation}\label{toto}
\hat{r}_h(x) = \frac{1}{n}\sum_{i=1}^n Y_{i} K_h(x - X_{i}) g^{-1}(X_{i}),
\end{equation}
with $K_h=K(\cdot/h)$ and it is shown that the estimator \eqref{toto} converges at the optimal rate \(n^{-\frac{\beta}{2\beta + 1}}\) over Hölder classes of regularity \(\beta > 0\) extending the result of \citet{lepski2014adaptive} that considered i.i.d. data.   \citet{mpadrino2024adaptiregre} constructed a GL estimator that satisfies an oracle-type inequality and is adaptive, converging at the rate  \(\left(\frac{n}{\log n}\right)^{-\frac{\beta}{2\beta + 1}}\) over Hölder classes of regularity  \(\beta > 0\). This rate is nearly optimal, up to the additional logarithmic factor. This factor is classical in the adaptive case and also appears, for example, in \citet{lepski2014adaptive}.

In this article, we extend the result of \citet{mpadrino2024adaptiregre} studying the estimation of the regression function $r$ at a fixed point $x\in \R$, in the case when the explanatory variable's density $g$ is unknown. Using that the regression function  \(r(\cdot) \) can be written as  \(r(\cdot)  = \frac{m(\cdot)}{g(\cdot)}\),  where for $x\in\R$ \(
m(x) = \int_{\mathbb{R}} y f_{X,Y}(x,y) dy\) with $f_{X,Y}$ is the joint density of the $(X_i,Y_i)$, we use a Nadaraya–Watson–type estimator $\hat{r}_{h/h^*}= \frac{\hat{m}_{h}}{\bar{g}_{h^*}}$
where for  \(h ,h^*> 0\)
$$\hat{m}_h(x)=\frac{1}{n}\sum_{i=1}^n K_{h}\left(x - X_i\right) Y_i$$ and $$\bar{g}_{h}(x) = \max\left\{ \hat{g}_{h}(x), \frac{1}{\log n} \right\},\quad\hat{g}_{h^*}= \frac{1}{n}\sum_{i=1}^n K_{h^*}\left(x - X_i\right). $$
The $1/\log n $ term is only here to prevent the denominator from vanishing in $r_{h/h^*}$. Such an approach has also been considered in \cite{comte}, but the properties of the estimators are derived in $L_2$ risk and i.i.d data.

Following the strategy outlined in \citet{mpadrino2024adaptiregre}, we show that  \(\hat{m}_h(\cdot)\) with $h=n^{-1/(2\beta_1+1)}$ converges at the optimal rate  \(n^{-\frac{\beta_1}{2\beta_1 + 1}}\) when $m$ belongs to H\"older classes of regularity \(\beta_1 > 0\). Then, selecting a bandwidth $\hat{h}$ by the GL method,  we prove that  \(\hat{m}_{\hat{h}}(\cdot)\) satisfies an oracle inequality and is adaptive, with a nearly optimal convergence rate  \(\left(\frac{n}{\log n}\right)^{-\frac{\beta_1}{2\beta_1 + 1}}\) over Hölder classes of regularity  \(\beta_1 > 0\). Similar minimax and adaptive results are satisfied for the estimator $\hat{g}_{h^*}(\cdot)
$ (see \citet{bertin2017pointwise}) with minimax rate \(n^{-\frac{\beta_2}{2\beta_2 + 1}}\) over Hölder classes of regularity \(\beta_2 > 0\).
 Using the Bernstein inequality for dependent data proposed in \citet{doukhan2007probability}, the following inequality is established:
\begin{equation} \label{IntroDes}
\E\left[\left(r(x) - \hat{r}_{h/h^*}(x)\right)^2\right] \leq F_1 \E\left[\left(m(x) - \hat{m}_{h}(x)\right)^2\right] + F_2 \E\left[\left(g(x) - \hat{g}_{h^*}(x)\right)^2\right] + \frac{F_3}{n^4},
\end{equation}
where \(F_1\), \(F_2\) and \(F_3\)  are constants. This inequality allows one to establish that the optimal convergence rate of the estimator  \(\hat{r}_{h/h^*}(\cdot)\) is \(n^{-\frac{\beta_{1} \wedge \beta_{2}}{2(\beta_{1} \wedge \beta_{2}) + 1}}\), where \(\beta_{1} \wedge \beta_{2}\) denotes the minimum between \(\beta_1\) and \(\beta_2\). Then, considering the family of estimators \(\{\hat{r}_{\hat{h}/\hat{h}^*}(\cdot), \; h \in \mathcal{H} \:\text{y}\: h^* \in \mathcal{H}^*\}\) for two families of bandwidths  \(\mathcal{H}\) and \(\mathcal{H}^*\), two bandwidths \(\hat{h}\) and \(\hat{h}^*\) are selected using the GL method to estimate  \(m(\cdot)\) and \(g(\cdot)\),  respectively. Subsequently, it is shown that the estimator  \(\hat{r}_{\hat{h}/\hat{h}^*}(\cdot)\) converges at a nearly optimal rate \(\left(\frac{n}{\log n}\right)^{-\frac{\beta_{1} \wedge \beta_{2}}{2(\beta_{1} \wedge \beta_{2}) + 1}}\).

Our contribution is the following: we obtain adaptive convergence rate for pointwise risks in a broad range of Hölder spaces, in the context of weak dependence, when the density of the explanatory variable is unknown. This extends the result presented in \citet{mpadrino2024adaptiregre}. Our estimator of $r$ is adaptive with respect to the smoothness of $m$ and $g$, in the sense that the two bandwidths are chosen from the data without knowledge of $\beta_1,\beta_2$. 

Moreover, our data-driven procedure depends only on explicit quantities, which means that it can be implemented directly in practice. As a consequence, we develop a new method for selecting a precise local bandwidth for kernel estimators. In particular, our method requires the calibration of two parameters,  \(\gamma\) and \(\gamma^*\). This calibration is carried out using unsupervised statistical learning techniques. A simulation study is also conducted to illustrate the performance of our method.

The remainder of this article is organized as follows: Section \ref{CT4Modelo} is devoted to presenting our model and the assumptions on the process \(\mathbb{X}\), including an example of a process that satisfies the stated hypotheses. The construction of our estimation procedure is developed in Section \ref{CT4EstPro}. The main results of the article are presented in Section \ref{CT4resultados}, which consists of three subsections: in \ref{CT4mchugchu}, we establish the bias, variance, consistency, and convergence rate of the estimator \(\hat{m}_{h}(\cdot)\), and we also state these results for the estimator \(\hat{g}_{h^*}(\cdot)\);  in \ref{CT4SubSecDOraculo}, we show that the kernel estimator  \(\hat{m}_{\hat{h}}(\cdot)\)  at the bandwidth selected by the GL method satisfies an oracle-type inequality, which allows us to conclude that this estimator is adaptive; and in \ref{CT4Resrhh}, we provide the inequality relation \eqref{IntroDes} between the MSE of \(\hat{r}_{h/h^*}(\cdot)\) and the mean squared errors of \(\hat{m}_{h}(\cdot)\) and \(\hat{g}_{h^*}(\cdot)\) (Proposition \ref{CT4DesConsVel}), which allows us to establish the consistency and convergence rate of the estimator \(\hat{r}_{h/h^*}(\cdot)\). Furthermore, we establish that the estimator \(\hat{r}_{\hat{h}/\hat{h}^*}(\cdot)\) at the bandwidths selected by the GL method is adaptive (Theorem \ref{CT4adapHatr}). 

A simulation study is carried out in Section \ref{CT4simulacion} to explain the technique used for calibrating the method and to illustrate its performance. 

\section{Model}\label{CT4Modelo}

We observe \(\{(X_1,Y_1),\ldots,(X_n,Y_n)\}\) identically distributed that satisfy  
\begin{equation}\label{modelt}
Y_i = r(X_i) + \varepsilon_i, \quad i=1,\ldots,n
\end{equation}
where \(n \ge 2\), the \(X_i\) are identically distributed with unknown density function \(g\), the \(\varepsilon_i\) are independent and identically distributed with distribution normal \(N(0,\sigma^2)\) with \(\sigma^2 > 0\), the sequence \(\varepsilon_i\) is independent of the process \(X_i\) , and the function \(r\) is the regression function \(r(x) = \mathbb{E}(Y \mid X=x)\).

Our goal is to estimate the function \( r \) at the point \( x \in \mathbb{R} \) using the observed sample \(\{(X_1,Y_1),\ldots,(X_n,Y_n)\}\) in a dependence approach.  
We will use the fact that the function \( r \) can be decomposed as
\begin{equation}\label{CT4fr}
r(x) = \mathbb{E}[Y \mid X = x] = \int_{\mathbb{R}} y f(y \mid x) dy = \int_{\mathbb{R}} y \frac{f_{X,Y}(x,y)}{g(x)} dy = \frac{m(x)}{g(x)},
\end{equation}
where \( f(\cdot \mid \cdot) \) is the conditional density function of \( Y \) given \( X \), which satisfies \( f(y \mid x) = \frac{f_{X,Y}(x,y)}{g(x)} \), with \( f_{X,Y} \) the joint density of \((X_i,Y_i)\), and
\begin{equation}\label{CT4eq:m}
m(x) = \int_{\mathbb{R}} y f_{X,Y}(x,y) dy.
\end{equation}
Given the expression \eqref{CT4fr}, as an estimation strategy, we will seek estimators of the functions \( m \) and \( g \) separately and take the quotient of both estimators as an estimator of \( r \).

The quality of the estimation of an estimator \(\hat{r}\) at a point \( x \in \mathbb{R} \) is measured using the mean squared error:
$$R(r,\hat{r},x)=\E \left(r(x)-\hat{r}(x)\right)^2.$$

On one hand, the following hypotheses are assumed, corresponding to bounds on the densities of the variables \(X_i\) in the neighborhood of \(x\): $$B(x)=[x-2/(\log n)^2, x+2/(\log n)^2].$$

\begin{itemize}
     \item[$(H_1)$] The density \( g \) of the \( X_i \) satisfies
      $$g_{\inf}\leq g(u) \leq g_{\sup}, \quad\forall u\in B(x)$$
where \( g_{\inf} \) and \( g_{\sup} \) are positive constants.
      \item[$(H_2)$] \(\forall (i,j)\), the joint density \(g_{i,j}\) of \((X_i, X_j)\) satisfies $$|g_{i,j}(u,v)|\leq Q, \quad\forall u,v\in B(x), $$
      where \(Q\) is a positive constant.
\end{itemize}

On the other hand, we assume a weak dependence structure on the variables $X_i$. 
More precisely, for $u$ and $v$ positive integers, we denote $i_{1:u}=(i_1,\ldots,i_u)\in\mathbb{Z}^u$ and $j_{1:v}=(j_1,\ldots,j_v)\in\mathbb{Z}^v$, then the random vectors $X_{i_{1:u}}=(X_{i_1},X_{i_2},\cdots, X_{i_u})$ and $X_{j_{1:v}} =(X_{j_1}, X_{j_2},\cdots, X_{j_v})$ are defined with values in $\mathbb{R}^u$ and $\mathbb{R}^v$, respectively. The function $q:\mathbb{Z}^u\times\mathbb{Z}^v \rightarrow \mathbb{Z}$ is defined by $q(i_{1:u},j_{1:v})=\min (j_{1:v}) - \max (i_{1:u})$ and $\Lambda_u$ is the class of functions $G_u:\mathbb{R}^u \rightarrow \mathbb{R}$ such that $\|G_u\|_{\infty} = \sup _{x\in\mathbb{R}^u} |G_u(x)| < \infty$. For a stochastic process $\mathbb{X}=(X_i)_{i\in\mathbb{Z}}$, the correlation coefficient $\alpha(\mathbb{X})=\left(\alpha_k(\mathbb{X})\right)_{k\in\mathbb{N}}$ is defined by,
$$\alpha_k (\mathbb{X})=\sup_{u,v\in\mathbb{N}} \sup_{(i_{1:u},j_{1:v})\in\mathbb{Z}^u\times\mathbb{Z}^v q(i_{1:u},j_{1;v})\geq k} \sup_{G_u \in \Lambda_u} \sup_{G_v \in \Lambda_v} \frac{|cov (G_u (X_{i_{1:u}}),G_v (X_{j_{1:v}} ))|}{\Psi (u,v,G_u , G_v )}$$
with $\Psi(u,v,G_u ,G_v )=4\|G_u\|_{\infty}\|G_v\|_{\infty}$.\\

We assume that the $\alpha$-mixing coefficient of the variables $X_i$ satisfies the following hypothesis.
\begin{itemize}
    \item[$(H_3)$] There exists $a \in ]0,1[$, such that
    $$\alpha_k (\mathbb{X})\leq a^k, \quad \forall k\in\mathbb{N}.$$
\end{itemize}

In the following, we will denote by $\mathcal{L}$ the set of processes $\mathbb{X}$ that satisfy $(H_1)$, $(H_2)$, and $(H_3).$\\


Next, two examples of processes that satisfy these hypotheses are given.
\begin{example}\label{CT4EjempAR1} 
    We consider $\mathbb{Z}=\{Z_t\}_{t\geq 1}$, an autoregressive process of order 1, defined by $Z_t =\phi Z_{t-1}+\rho \xi_t$, with $\xi_t \sim N(0,1)$, $|\phi|<1$, and $\rho >0$. By recurrence, $Z_{t}=\phi^{n}Z_{t-n}+\rho\sum_{i=0}^{n-1}\phi^{i}\xi_{t-i}$ and taking $n\rightarrow \infty$, we have $Z_{t}=\rho\sum_{i=0}^{\infty}\phi^{i}\xi_{t-i}$, resulting in a centered and stationary Gaussian process with covariance function $\gamma_{Z}(k)=cov(Z_{t},Z_{t+k})=\frac{\rho^2}{1-\phi^2}\phi^k$ for $k\geq 0$. It can be shown that the process $\mathbb{Z}$ is weakly dependent, thus satisfying Hypothesis $(H_{3})$, see \citet{doukhan1994mixing}. 
\end{example}

\begin{example}\label{CT4ProcesoXenL}
    Denoting $\phi_{\mu,\sigma^2}(x)$ as the density function $N(\mu,\sigma^{2})$ and $\Phi_{\mu,\sigma^2}(x)$ as the cumulative distribution function $N(\mu,\sigma^{2})$. For $c\in\mathbb{R}^+$, we consider the density and cumulative distribution functions of  truncated normal distribution in $[-c,c]$ given by $g(x)=\frac{\phi_{0,1}(x)}{p}\mathbf{1}_{[-c,c]}(x)$, with $p=\Phi_{0,1}(c)-\Phi_{0,1}(-c)$, for each $x\in\mathbb{R}$ and
    $$G(x)=\left\{ \begin{array}{lcc}
         0 & \text{if} & x<-c  \\
         \\ \frac{1}{p}\left(\Phi_{0,1}(x)-\Phi_{0,1}(-c) \right)& \text{if} & -c\leq x < c \\
         \\ 1 & \text{if} & x\geq c
    \end{array}
    \right.$$
    respectively. Furthermore, the inverse function of $G$ is $G^{-1}(u)=\Phi_{0,1}^{-1}\left( pu+\Phi_{0,1}(c) \right)$ for each $u\in[0,1]$.

    It can be shown that the process $\mathbb{X}=\{X_t\}_{t\geq 1}$ defined by $X_{t}=\left(G^{-1}\circ\Phi_{0,\frac{\rho^2}{1-\phi^2}}\right)(Z_{t})$ satisfies that $\mathbb{X}\in\mathcal{L}$, where the process $\mathbb{Z}=\{Z_{t}\}_{t\geq 1}$ with each 
    $Z_t$ having distribution $\Phi_{0,\frac{\rho^2}{1-\phi^2}}$ is the process from Example \ref{CT4EjempAR1}.\\
\end{example}

Finally, we assume the following hypothesis about the regression function.
\begin{itemize}
        \item[$(H_4)$] There exists a positive constant $r_{\sup}$ such that $$|r(u)|\leq r_{\sup}, \quad \forall u\in B(x).$$
\end{itemize}

Note that since $m(\cdot)=r(\cdot)g(\cdot)$, the hypotheses $(H_1)$ and $(H_4)$ imply that
$\forall u\in B(x)$
$$|m(u)|\leq r_{\sup}g_{\sup}.$$

\section{Statistical Procedure}
\label{CT4EstPro}

For the sample $\{(X_i,Y_i)\}_{i=1}^{n}$, we consider the kernel estimators of the density function $g$ and the joint density $f_{X,Y}$ given for $h>0$ by:

\begin{equation}\label{CT4eq:gh}
\hat{g}_{h}(x)=\frac{1}{n}\sum_{i=1}^nK_{h}\left(x-X_i\right)
\end{equation}
and
\begin{equation*}
\hat{f}_{X,Y}^{(h)}(x,y)=\frac{1}{n}\sum_{i=1}^nK_{h}\left(x-X_i\right)K_{h}\left(y-Y_i\right)
\end{equation*}
where $K_{h}(\cdot)=\frac{1}{h}K(\frac{\cdot}{h})$ and $K:\mathbb{R}\rightarrow \mathbb{R}$ is a kernel function satisfying $\int K(u)du =1$ and the following hypothesis.
\begin{itemize}
    \item[$(H_5)$] $K$ has support $[-1,1]$ and $\|K\|_{\infty}<+\infty$.
    \end{itemize}
This hypothesis implies that $\|K\|_{1}<+\infty$ and $\|K\|_{2}<+\infty$.

By substituting in the equation \eqref{CT4eq:m}, the joint density with $\hat{f}_{X,Y}^{(h)}(x,y)$ (and assuming $\int uK(u)du=0$), the estimator of $m$ is obtained.
\begin{equation} \label{CT4m_est}
\hat{m}_{h}(x)=\frac{1}{n}\sum_{i=1}^nK_{h}\left(x-X_i\right)Y_i.
\end{equation}

Our estimation strategy consists of taking the quotient of the estimators \eqref{CT4eq:gh} and \eqref{CT4m_est}. In order to avoid $\hat{g}_{h}(x)$ being zero, we define an alternative estimator of $g$ defined by
\begin{equation} \label{CT4Barg_h}
    \bar{g}_{h}(x) = \max\left\{ \hat{g}_{h}(x), \frac{1}{\log n} \right\}.
\end{equation}
Finally, we consider the estimator of $r$.
\begin{equation} \label{CT4Hatr_hh*}
\hat{r}_{h/h^*}=\frac{\hat{m}_{h}}{\bar{g}_{h^*}}
\end{equation}
which depends on two bandwidths $h, h^* > 0$. Such an approach has also been considered in \cite{comte}. Observe that if $h = h^*$, this estimator is related to the Nadaraya-Watson estimator (\cite{NW}).  \\




To select the bandwidths $h$ and $h^*$ of the estimator, we will use the Goldenshluger-Lepski (GL) method and select $h$ and $h^*$ separately for the estimator of $m$ and the estimator of $g$. \\

\textbf{Bandwidth selection procedure to estimate $m$:} We consider the family of bandwidths
$$\mathcal{H}=\{e^{-i}\}_{i=0}^{M}\cap[h_{\min},h_{\max}]$$ 
with $h_{\min} = \frac{(\log n)^8}{n}$, $h_{\max} = \frac{1}{(\log n)^2}$, and $M = \left[\log\left(\frac{n}{(\log n)^8}\right)\right]$. An auxiliary smoothing estimator is defined for $h, h' \in \mathcal{H}$
\begin{equation}\label{CT4nwss}
\hat{m}_{h,h'}(x)=\frac{1}{n}\sum_{i=1}^nK_{h}\ast K_{h'}\left(x-X_i\right)Y_i.
\end{equation}\label{CExpresionA}
We consider for $h \in \mathcal{H}$,
\begin{equation}
    A(h,x)= \max_{h'\in\mathcal{H}}\{|\hat{m}_{h,h'}(x)-\hat{m}_{h'}(x)|-V(h')\}_+
\end{equation}
where $\{ y \}_{+} = \max(0, y)$ denotes the non-negative part of $y \in \mathbb{R}$ and $V(h)$ is given by,

\begin{equation}\label{CExpresionV}
    V(h)=\sqrt{ 2\gamma A_{1}}( \|K\|_1 + 1) ( 1 + \delta_n ) \frac{(\log n)^\frac{1}{2}}{(nh)^\frac{1}{2}}
\end{equation}
with  $\gamma>2$,  $A_{1} = ((r_{\sup})^{2}+\sigma^2)g_{\sup}\|K\|_2^2$ and $\delta_n =(\log n)^{-\frac{1}{5}}$.\\

The procedure consists of selecting, based on the data, a bandwidth $\hat{h}$ from the family $\mathcal{H}$, given by
$$
\hat{h} = \mbox{argmin}_{h\in\mathcal{H}}\{A(h,x)+V(h)\}
$$
We show that the resulting estimator, $\hat{m}_{\hat{h}}$, satisfies an oracle inequality that allows us to demonstrate that it is adaptive. \\

\textbf{Bandwidth selection procedure to estimate $g$:} Here we will use the method proposed in \cite{bertin2017pointwise}. Analogously, a family of bandwidths is defined
$$\mathcal{H}^{*}=\{2^{-k}:k\in\mathbb{N}\} \cap [h_{\min}^{*},h_{\max}^{*}]$$
where $h_{\min}^{*} = n^{-1}\exp{(\sqrt{\log n} - 1)}$ and $h_{\max}^{*} = \exp{(-\sqrt{\log n})}$.\\

We consider for $h,h' \in \mathcal{H}^{*}$
\begin{equation}\label{CExpresionA*}
    A^{*}(h,x) = \max_{h'\in\mathcal{H}^*} \left\{|\hat{g}_{h\vee h'}(x)-\hat{g}_{h'}(x)|-V^{*}(h,h') \right\}_{+}
\end{equation}

\begin{equation}\label{CExpresionV*}
    V^{*}(h) = \sqrt{2\gamma^{*}|\log h|\left( \hat{J}_{n}(h)+\frac{\delta^*_{n}}{nh} \right)}, \quad V^{*}(h,h') = V^{*}(h')+V^{*}(h \vee h')
\end{equation}
where $\delta^*_{n}=(\log n)^{-\frac{1}{2}}$, $h\vee h' = \max\{h,h'\}$ and
$\hat{J}(h) = \frac{1}{n^2}\sum_{i=1}^{n} K^2_{h}(x-X_{i})$.

The procedure consists of selecting, based on the data, a bandwidth $\hat{h}^*$ from the family $\mathcal{H}^*$, given by

$$
\hat{h}^{*} = \mbox{argmin}_{h \in \mathcal{H}^*}\{A^{*}(h,x) + V^{*}(h)\}
$$
The resulting estimator, $\hat{g}_{\hat{h}^*}$, satisfies an oracle inequality and is adaptive (See Theorems 1 and 2 in \citet{bertin2017pointwise}).\\

\section{Results} \label{CT4resultados}
\vspace{0.5cm}
\subsection{Results of the estimators $\hat{m}_{h}$ and $\hat{g}_{h}$.}
\label{CT4mchugchu}

In this section, we study the bias, variance, and convergence rate of the estimator $\hat{m}_{h}$, and we state these same results for the density estimator $\hat{g}_{h}$. Let us first define Hölder classes and kernels of order $\mathit{m}_0$.

\begin{definition}[Hölder Class] Let $\beta > 0$, $L > 0$. The Hölder class $\Sigma(\beta, L)$ is defined as the set of all functions $f:\mathbb{R} \rightarrow \mathbb{R}$ such that the derivative $f^{\mathit{(l)}}$, $\mathit{l} = \left \lfloor \beta \right \rfloor$ exists and 
\begin{equation*}
\left|f^{(\mathit{l})}(x) - f^{(\mathit{l})}(y)\right| \leqslant L\left|x-y\right|^{\beta - \mathit{l}}, \forall x, y \in \mathbb{R}
\end{equation*}
where $\left \lfloor \beta \right \rfloor = \max\left\{ n \in \mathbb{N}, n < \beta \right\}$.
\end{definition}

\begin{definition}[Kernel of order $m_0$] Let $m_0 \in \mathbb{N}$. A function $K:\mathbb{R} \rightarrow \mathbb{R}$ is said to be a kernel of order $m_0$ if the functions $u \rightarrow u^{j}K(u)$, for $j=0,1,...,m_0$, satisfy 
\begin{equation*}
\int_\mathbb{R} K(u) du = 1, \quad \int_\mathbb{R} u^{j} K(u) du = 0, \quad j=1,...,m_0.
\end{equation*}
\end{definition}
Next, in the following two propositions, we bound the bias and variance of the estimator $\hat{m}_h$ (proofs in Appendix \ref{DemoSesgoHatm} and \ref{DemoVarHatm}).

\begin{proposition}\label{CT4pro:sesgo} 
    Let $\beta_1 > 0$ and $L_1 > 0$. We assume that $K$ is of order $m_0$, where $m_0 \geq \mathit{l}_1 = \left \lfloor \beta_1 \right \rfloor$, and that it satisfies $\int_\mathbb{R} |u|^{\beta_1} |K(u)| du < \infty$. Then, if $m \in \Sigma(\beta_1, L_1)$, we have
$$\left|E\left[\hat{m}_{h}(x)\right] - m(x)\right| = \left|K_h \ast m (x) - m(x)\right| \leq h^{\beta_1} A_{0},$$
where $A_0 = \frac{L_1}{\mathit{l}_1!} \int_{\mathbb{R}} |u|^{\beta_1} |K(u)| du$. 
\end{proposition}
\begin{proposition}\label{CT4pro:var}
 Under the hypotheses $(H_1)$, $(H_2)$, $(H_3)$, $(H_4)$, and $(H_5)$, it holds that for $h \in (0, h_{max})$
  $$\mathbb{E}\left[(\hat{m}_h (x)-\mathbb{E}[\hat{m}_h (x)])^2\right] \leq \frac{A_{1}}{nh}+A_{2}\frac{(\log n)^{-\frac{1}{2}}}{nh},$$
 for each $n \geq 4$, where the constants are defined by $A_{1}=((r_{\sup})^2 +\sigma^2 )g_{\sup}\|K\|_2^2$, $A_{2} = \frac{2A_3}{|\log a|} + 20A_4$, $A_3 = (r_{\sup})^2 \|K\|_1^2 \left(Q+(g_{\sup})^{2}\right)$, and $A_4 =\frac{4}{1-a}(r_{\sup})^2 \|K\|_{\infty}^2$.
\end{proposition}

Combining Propositions \ref{CT4pro:sesgo} and \ref{CT4pro:var} we derive the rate of convergence of the estimator $\hat{m}_h$ in the following theorem.

\begin{theorem}\label{CT4veloconv}
Under the hypotheses of Propositions \ref{CT4pro:sesgo} and \ref{CT4pro:var}, if $m \in \Sigma(\beta_1, L_1)$, then the estimator $\hat{m}_{\tilde{h}}$, with 
$\tilde{h} = n^{-\frac{1}{2\beta_1 + 1}}$, satisfies
$$\E\left[\left(\hat{m}_{\tilde{h}}(x) - m(x)\right)^{2}\right] \leq C n^{-\frac{2\beta_1}{2\beta_1 + 1}},$$
where $C$ is a positive constant that depends on the kernel $K$, on $\beta_1$, $L_1$, $r_{sup}$, $Q$, $g_{sup}$, $\sigma$, and $a$.
\end{theorem}


Following  \citet{bertin2017pointwise} similar results are satisfied for the estimator \(\hat{g}_{h}\) that we gather in the following proposition.

\begin{proposition}[\citet{bertin2017pointwise}]\label{resume} 
   Let $\beta_2>0$ and $L_2>0$. Let $g\in\Sigma(\beta_2,L_2)$. We assume that $K$ is of order $m_0$, where $m_0\geq \mathit{l}_2=\left \lfloor \beta_2 \right \rfloor$, and that it satisfies $\int_\mathbb{R}|u|^{\beta_2}|K(u)|du<\infty$. Assume also that assumptions $(H_1)$, $(H_2)$, $(H_3)$, and $(H_5)$ are satisfied. Then we have for \( h \in (0, h_{max}^*) \)
   
\begin{equation}\label{CT4pro:sesgoHatgh}\left|\E\left[\hat{g}_{h}(x)\right]-g(x)\right|\leq h^{\beta_2}\mathfrak{A}_{0},
\end{equation}
\begin{equation}\label{CT4pro:varHatgh}
\mathbb{E}\left[(\hat{g}_h (x)-\mathbb{E}[\hat{g}_h (x)])^2\right] \leq \frac{\mathfrak{A}_1}{nh}+\mathfrak{A}_2\frac{(\log n)^{-\frac{1}{2}}}{nh},
\end{equation}
where $\mathfrak{A}_0 =\frac{L_2}{\mathit{l}_2 !}\int_{\mathbb{R}}|u|^{\beta_2}|K(u)|du$, \( \mathfrak{A}_1 = g_{\sup}\|K\|_2^2 \) and \( \mathfrak{A}_2 = \frac{2\mathfrak{A}_3}{|\log a|} + 20\mathfrak{A}_4 \) with \( \mathfrak{A}_3 = \|K\|_1^2 \left(Q+(g_{\sup})^{2}\right) \) and \( \mathfrak{A}_4 = \frac{4}{1-a}\|K\|_{\infty}^2 \).

Combining \eqref{CT4pro:sesgoHatgh} and \eqref{CT4pro:varHatgh}, we have
\begin{equation}\label{CT4veloconvHatgh}\E\left[\left(\hat{g}_{\tilde{h}^{*}}(x)-g(x)\right)^{2}\right]\leq \mathfrak{C} n^{-\frac{2\beta_2}{2\beta_{2}+1}},
\end{equation}
where \(\tilde{h}^{*}=n^{-\frac{1}{2\beta_2 +1}}\) and \(\mathfrak{C}\) is a positive constant that depends on the kernel \(K\), on \(\beta_2\), \(L_2\), \(Q\), \(g_{sup}\), and \(a\).
\end{proposition}
 
 \subsection{Oracle Inequality and Adaptivity of the Estimator $\hat{m}_{\hat{h}}$}\label{CT4SubSecDOraculo}

Now we consider the estimator $\hat{m}_{\hat{h}}$ where $\hat{h}$ has been selected by the GL method defined in Section~\ref{CT4EstPro}. In the two following theorems, we obtain that it satisfies an oracle inequality and that it converges at the adaptive rate of convergence $\left(\frac{n}{\log n}\right)^{-\frac{\beta_1}{2\beta_{1} +1}}$ on H\"older classes of regularity $\beta_1$.
\begin{theorem} \textbf{(Oracle Inequality)}\label{CT4teo:D-Oraculo}
Under assumptions $(H_1)$, $(H_2)$, $(H_3)$, $(H_4)$ and $(H_5)$, the estimator $\hat{m}{\hat{h}}(x)$ satisfies the following inequality,
\begin{eqnarray} \label{CT4eq:oracle}
\nonumber
\left( \mathbb{E}\left(m(x)-\hat{m}_{\hat{h}}(x)\right)^2\right)^{\frac{1}{2}} &\leq& \min_{h\in\mathcal{H}} \left( A_6 C(h) + A_7 ( 1 + \delta_n ) \frac{(\log n)^{\frac{1}{2}}}{(nh)^{\frac{1}{2}}} \right) + A_8\frac{(\log n)^\frac{1}{2}}{n^\frac{1}{2}}\
\end{eqnarray}
where $A_6 =1+2\|K\|_1$, $A_7 =\sqrt{A{1} +A_2} +2\sqrt{ 2\gamma A_{1}}( \|K\|_1 + 1)$ and $A_8$ is a positive constant depending on the kernel $K$, on $r_{sup}$, $Q$, $g_{sup}$, $\sigma$, $a$, and $\gamma$.
\end{theorem}
The proof of Theorem~\ref{CT4teo:D-Oraculo} can be found in Appendix \ref{DemoOracuHatm}. Now, using Theorem~\ref{CT4teo:D-Oraculo} and following the same proof as in Theorem 3 of Section 4.2 in \citet{mpadrino2024adaptiregre}, we deduce the following theorem for the estimator $\hat{m}_{\hat{h}}$.

\begin{theorem} \textbf{(Adaptivity)} \label{CT4AdapHatm} We assume that assumptions $(H_1)$, $(H_2)$, $(H_3)$, $(H_4)$ and $(H_5)$ hold. We also assume that $K$ is of order $m_0$. Then, for all $0<\beta_1\leq m_0$ and $L_1>0$, the estimator $\hat{m}_{\hat{h}}$ satisfies, for $m\in\Sigma(\beta_1,L_1)$,
\begin{equation*}
\left(\mathbb{E}\left(m(x)-\hat{m}_{\hat{h}}(x)\right)^2\right)^{\frac{1}{2}} \leq C^*\left(\frac{n}{\log n}\right)^{-\frac{\beta_1}{2\beta_{1} +1}}
\end{equation*}
where $C^*$ is a positive constant depending on the kernel $K$, on $\beta_1$, $L_1$, $r_{sup}$, $Q$, $g_{sup}$, $\sigma^2$, $a$, and $\gamma$.
\end{theorem}

\begin{remark}
As classically in oracle inequalities, in Theorem \ref{CT4teo:D-Oraculo}, we can see that the estimator $\hat{m}_{\hat{h}}$ mimics the "oracle"; i.e., the best possible (but unknown) estimator $\hat{m}_{h}$ in the family $(\hat{m}_h)_{h\in\mathcal{H}}$, which minimizes the sum of bias $C(h)$ and standard deviation. \citet{comte} obtained a similar oracle inequality (Theorem 3.2) in an i.i.d. context and in $L_2$ risk.
The adaptive estimator $\hat m_{\hat h}$ attains the rate $(n/\log n)^{-\beta_1/(2\beta_1+1)}$ over Hölder classes of regularity $\beta_1$, without knowing $\beta_1$. This is minimax-optimal up to a logarithmic factor, which is the usual price of adaptation for pointwise risk (see \citet{lepski1990problem}). To our knowledge, this pointwise adaptive result is new in the setting of weakly dependent data with an unknown design density.
\end{remark}

\subsection{Results for the estimator $\hat{r}_{h/h^*}$} \label{CT4Resrhh}

We now study the properties of the estimator $\hat{r}_{h/h^*}$. This estimator satisfies
\begin{align}
r(x)-\hat{r}_{h/h^*}(x)=&\frac{m(x)}{g(x)}-\frac{\hat{m}_
{h}(x)}{\bar{g}_{h^*}(x)}=
\frac{m(x)}{\bar{g}_{h^*}(x)}-\frac{\hat{m}_{h}(x)}{\bar{g}_{h^*}(x)}+\frac{m(x)}{g(x)}-\frac{m(x)}{\bar{g}_{h^*}(x)}\nonumber\\
=&
\frac{1}{\bar{g}_{h^*}(x)}\left(m(x)-\hat{m}_{h}(x)\right)+\frac{m(x)}{g(x)}\left(1-\frac{g(x)}{\bar{g}_{h^*}(x)}\right)\nonumber\\
=&\frac{1}{\bar{g}_{h^*}(x)}(m(x)-\hat{m}_{h}(x))+\frac{r(x)}{\bar{g}_{h^*}(x)}(\bar{g}_{h^*}(x)-g(x)).\label{CT4decomp}
\end{align}

From formula \eqref{CT4decomp}, we see that the properties of $\hat{r}_{h/h^*}$ can be deduced from the properties of the estimators of $m$ and $g$. We have the following proposition proved in Appendix \ref{DemoPropo5}.

 \begin{proposition} \label{CT4DesConsVel}
 We assume that hypotheses $(H_1)$, $(H_2)$, $(H_3)$, $(H_4)$, and $(H_5)$ hold. Let $\beta_2>0$ and $L_2>0$. Assume that $g\in\Sigma(\beta_2,L_2)$, and that $K$ is of order $m_0\geq\lfloor\beta_2\rfloor${ and $n$ is sufficiently large}. For $h\in[h_{\min},h_{\max}]$ and $h^*\in [h^*_{\min},h^*_{\max}]$, we have
     $$\E\left[\left(r(x)-\hat{r}_{h/h^*}(x)\right)^2\right]\leq F_1\E\left[\left(m(x)-\hat{m}_{h}(x)\right)^2\right]+F_2\E\left[\left(g(x)-\hat{g}_{h^*}(x)\right)^2\right]+\frac{F_3}{n^4}$$
     where $F_1=\frac{8}{(g_{\inf})^2}$, $F_{2}=8\frac{(r_{\sup})^2}{(g_{\inf})^2}$ and $F_3$ is a positive constant that depends on the kernel $K$, on $\beta_2$, $L_2$, $r_{\sup}$, $Q$, $g_{\inf}$, $g_{\sup}$, and $a$.
\end{proposition}

Under the assumptions of Propositions \ref{CT4pro:sesgo}, \ref{CT4pro:var}, \ref{resume}, and \ref{CT4DesConsVel},  the estimator $\hat{r}_{h/h^*}$ is a consistent estimator of $r$, as soon as  $n\rightarrow \infty$, $h\rightarrow 0$, $h^* \rightarrow 0$, $nh \rightarrow \infty$ and $nh^* \rightarrow \infty$. Moreover we can derive its rate of convergence.

\begin{theorem} \textbf{(Rate of convergence)} \label{CT4VelConvHatr}
Under the assumptions of Theorem \ref{CT4veloconv} and Proposition \ref{resume}, if $m\in\Sigma(\beta_{1},L_{1})$ and $g\in\Sigma(\beta_{2},L_{2})$, then there exists a constant $\mathfrak{D}>0$ such that the estimator $\hat{r}_{\tilde{h}/\tilde{h}^*}$ with $\tilde{h}=n^{-\frac{1}{2\beta_{1}+1}}$ and $\tilde{h}^{*}=n^{-\frac{1}{2\beta_2 +1}}$ satisfies

     $$\E\left[\left(\hat{r}_{\tilde{h}/\tilde{h}^*}(x)-r(x)\right)^2\right] \leq \mathfrak{D}n^{-\frac{2(\beta_{1} \wedge \beta_{2})}{2(\beta_{1} \wedge \beta_{2})+1}}.$$
where $\mathfrak{D}$ is a positive constant depending on the kernel $K$, on $\beta_1$, $L_1$, $\beta_2$, $L_2$, $r_{\sup}$, $g_{\inf}$, $g_{\sup}$, $Q$, $\sigma$, $a$, and ${\alpha_{k}(\mathbb{X}) }_{k\geq 0}$.
\end{theorem}
Proposition \ref{resume} and \ref{CT4DesConsVel}, and   Theorem \ref{CT4veloconv} imply that
$$\E\left[\left(\hat{r}_{\tilde{h}/\tilde{h}^*}(x)-r(x)\right)^2\right] \leq F_1C n^{-\frac{2\beta_1}{2\beta_{1}+1}}+F_2\mathfrak{C} n^{-\frac{2\beta_2}{2\beta_{2}+1}}+\frac{F_3}{n^4},$$
which yields the result of the theorem.

We now study the estimator $\hat{r}_{\hat{h}/\hat{h}^*}=\frac{\hat{m}_{\hat{h}}}{\bar{g}_{\hat{h}^*}}$, where $\hat{h}$ y $\hat{h}^*$ are defined in Section \ref{CT4EstPro}. We obtain the following theorem (see proof in Appendix \ref{DemoCT4adapHatr}).

\begin{theorem} \textbf{(Adaptivity)} \label{CT4adapHatr}
We assume that hypotheses $(H_1)$, $(H_2)$, $(H_3)$, $(H_4)$ and $(H_5)$ hold. We assume that $K$ is of order $m_0$. Then for all $0<\beta_1\leq m_0$, $L_1>0$, $0<\beta_2\leq m_0$, and $L_2>0$, the estimator $\hat{r}_{\hat{h}/\hat{h}^*}$ satisfies, for $m\in\Sigma(\beta_1,L_1)$ and $g\in\Sigma(\beta_{2},L_{2})$,
    $$\left( \E\left[\left(\hat{r}_{\hat{h}/\hat{h}^*}(x)-r(x)\right)^2\right] \right)^{\frac{1}{2}} \leq \mathfrak{B}\left( \frac{n}{\log n} \right)^{-\frac{(\beta_{1} \wedge \beta_{2})}{2(\beta_{1} \wedge \beta_{2})+1}}$$
    where $\mathfrak{B}$ is a positive constant depending on the kernel $K$, on $\beta_1$, $L_1$, $\beta_2$, $L_2$, $r_{sup}$, $Q$, $g_{inf}$, $g_{sup}$, $\sigma^2$, and $a$.
\end{theorem}
\begin{remark}
We find that the estimator $\hat{r}_{\hat{h}/\hat{h}^*}$ converges at the rate $\left( \frac{n}{\log n} \right)^{-\frac{(\beta_{1} \wedge \beta_{2})}{2(\beta_{1} \wedge \beta_{2})+1}}$ and that it is adaptive with respect to the smoothness of $m$ and $g$, in the sense that the two bandwidths are chosen from the data without knowledge of $\beta_1$ and $\beta_2$. \citet{comte} obtain the same rate of convergence in $L_2$ risk and for i.i.d. data. If we assume that $r$ is $\beta$-H\"older with $\beta>0$, the rate in Theorem~\ref{CT4adapHatr} coincides with the optimal adaptive rate when $\beta=\beta_{1} \wedge \beta_{2}$.
\end{remark}

\section{Simulation Study} \label{CT4simulacion}

\vspace{0.5cm}

In this section, we illustrate the performance of the proposed procedure using simulated data. Global and local empirical errors are compared for different sample sizes. Finally, the calibration scheme used for the method is presented in Appendix~\ref{Calibra} in order not to interrupt the main exposition. The objective is to estimate the regression function when the explanatory variable is a weakly dependent process with unknown density.

\subsection{Simulation Scheme} \label{CAj}

We consider the regression model \eqref{modelt} with the regression function $r(x)=0.7x+2e^{-10x^2}$, which, when restricted to the interval $[-2,2]$, satisfies Assumption $(H_{4})$.
We consider that the variables $\varepsilon_i$ are independent and identically distributed with distribution $N(0,\sigma^2)$, with $\sigma\in \{0.1,0.5,1\}$. 

Let $n,q\in\mathbb{N}$. A sample $\{X_i\}_{i=1}^{n+2q}$ is generated from the process $\mathbb{X}$ introduced in Example~\ref{CT4ProcesoXenL}, where
$X_{i}=\left(G^{-1}\circ\Phi_{0,\frac{\rho^2}{1-\phi^2}}\right)(Z_{i})$
with $\rho=1$, $\phi=0.75$, and $c=2$.  The variables $X_i$ follow a truncated normal distribution with the density function
$g(x)=\frac{\phi_{0,1}(x)}{p}\mathbf{1}_{[-2,2]}(x)$, where
$p=\Phi_{0,1}(2)-\Phi_{0,1}(-2)$, for each $x\in\mathbb{R}$. As mentioned in Example~\ref{CT4ProcesoXenL}, the process satisfies $\mathbb{X}\in\mathcal{L}$, that is, assumptions $(H_{1})$, $(H_{2})$, and $(H_{3})$ hold.

\begin{figure}[htp!] 
	\centering 
	\begin{subfigure}
	\centering
	\includegraphics[width=7.66cm]{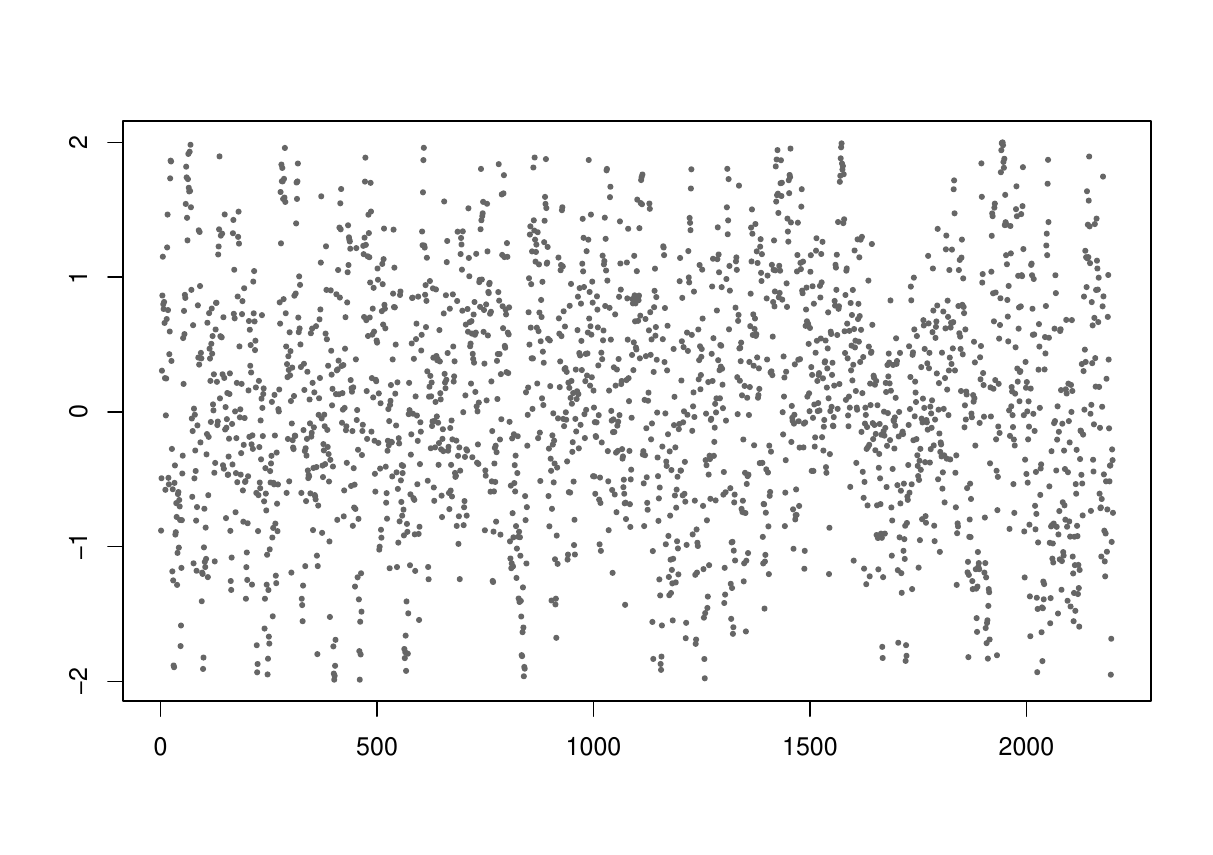}  
	\end{subfigure}
	\begin{subfigure}
	\centering
\includegraphics[width=7.66cm]{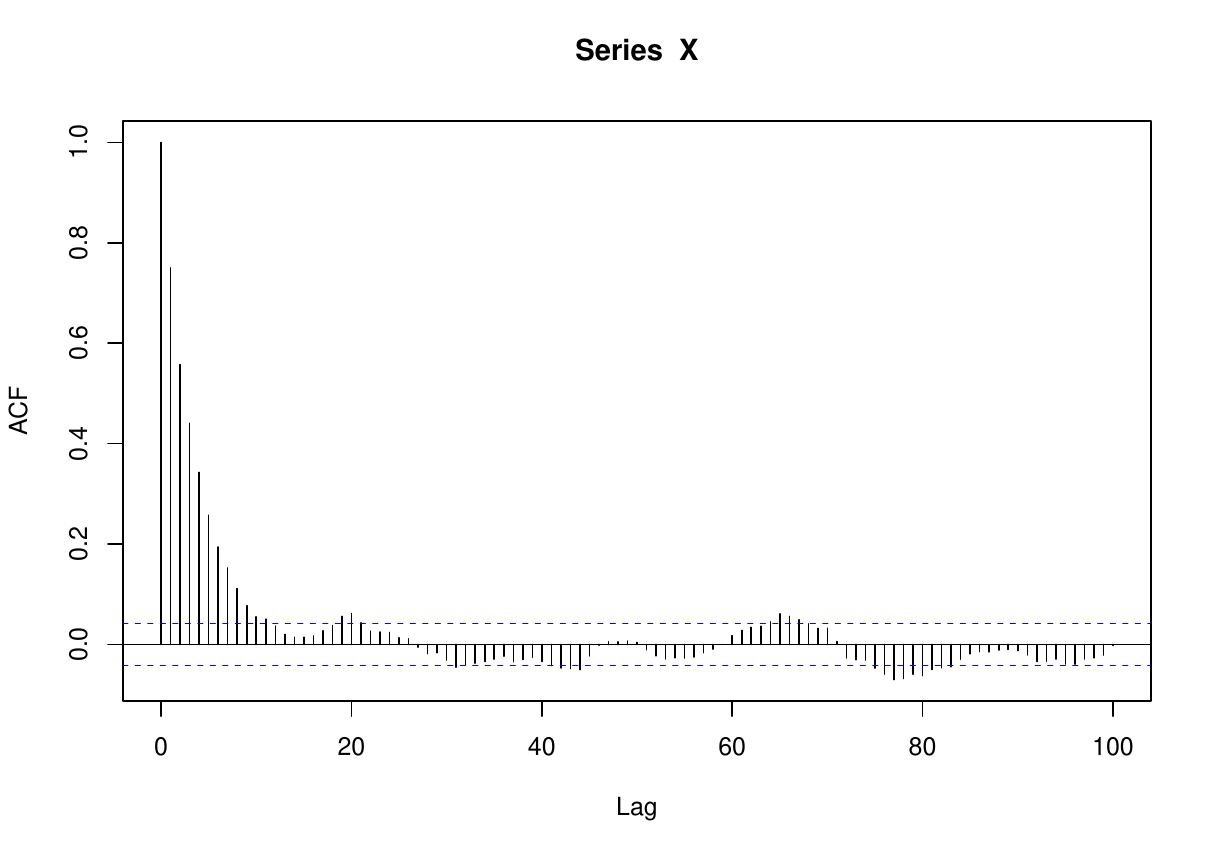}  
	\end{subfigure}
	\caption{Scatter plot of $\mathbb{X}$ and sample autocorrelation function $\gamma_{X}$.}\label{CXyCovaX2000}
\end{figure}

Figure~\ref{CXyCovaX2000} displays the scatter plot of the observations $\{X_i\}_{i=1}^{n+2q}$ for $n=2000$ and $q=100$ (left panel), together with the sample autocorrelation function $\gamma_{X}(k)=corr(X_{i},X_{i+k})$ for $k=1,\cdots,100$ (right panel). The exponential decay of the correlations can be observed, a characteristic feature of weakly dependent $\alpha$-mixing processes.\\

The first $n$ observations $\{(X_i,Y_i)\}_{i=1}^{n}$ constitute the sample used in Section~\ref{CEstGL} to estimate the regression function, while the last $q$ observations $\{(X_i,Y_i)\}_{i=n+q+1}^{n+2q}$ are used to calibrate the method by selecting suitable values of $\gamma^*$ and $\gamma$, as described in Appendix~\ref{Calibra}. The purpose of introducing a time gap of length $q$ between the estimation sample and the calibration sample is to ensure that the $\alpha_m(\mathbb{X})$-mixing coefficient is sufficiently small, thereby reducing the overfitting effect induced by the dependence between the two samples. Control of the value of the $\alpha_m(\mathbb{X})$-mixing coefficient is ensured through the exponential decay assumption $(H_{3})$.

\subsection{Adaptive Estimation of the Regression Function} \label{CEstGL}

The regression function $r$ is estimated on the equally spaced grid $\{x_j\}_{j=1}^{s}$ over the interval $[-1,1]$, where
$x_{j}=-1+2\left(\frac{j-1}{s-1}\right)$ for $j=1,\cdots,s$, and we set $s=21$. The estimator of the regression function $r$ at the point $x_j$ is
$\hat{r}_{h_j/h^*_j}(x_j )=\frac{\hat{m}_{h_j}(x_j)}{\hat{g}_{h^*_j}(x_j )}$
where for $h_{j},h^*_j\in\mathcal{H}$,
$$
\hat{m}_{h_j}(x_{j})=\frac{1}{n}\sum_{i=1}^nK_{h_j}\left(x_{j}-X_i\right)Y_i \quad\text{and}\quad \hat{g}_{h^*_{j}}(x_{j})=\frac{1}{n}\sum_{i=1}^{n}K_{h^*_j}(x_{j}-X_{i})
$$
and $\mathcal{H}$ denotes the collection of bandwidths $\mathcal{H}=\left\{e^{-0.1j} : j=0,\cdots,\frac{[\log(n)]^{\frac{2}{3}}}{0.1} \right\}$. The kernel
$K(x)=\frac{1}{\sqrt{2\pi}}e^{-\frac{x^2}{2}}$, $x\in\mathbb{R}$,
is employed. Although it is not compactly supported, it exhibits good practical performance. Moreover, it is well known that
$\|K\|_{1}=1$ and $\|K\|_{2}=\frac{1}{\sqrt{2\sqrt{\pi}}}$.\\
We aim at studying the performance of the estimators $\hat{r}_{\hat{h}_{j}/\hat{h}^*_j}(x_{j})$ for $j=1,\cdots,s$, where
the bandwidths $\hat{h}_{j}$ and $ \hat{h}^*_j$ are selected by the GL method described below. Note that since we adopt a pointwise framework, the bandwidths depend on the estimation point.

The bandwidth $\hat{h}^*_j$ is selected from $\mathcal{H}$ by taking
\begin{equation}\label{Ch*xi}
\hat{h}^*_{j}=\mbox{arg}\hspace{-0.05cm}\min_{h\in\mathcal{H}}\{A^{*}(h,x_{j})+V^{*}(h)\}
\end{equation}
for each $j=1,\cdots,s$, where
$A^{*}(h,x_{j})= \max_{h'\in\mathcal{H}}\{|\hat{g}_{h,h'}(x_{j})-\hat{g}_{h'}(x_{j})|-V^{*}(h')\}_+$,
with
$\hat{g}_{h,h'}(x_{j})=\frac{1}{n}\sum_{i=1}^nK_{h}\ast K_{h'}\left(x_{j}-X_{i}\right)$
and
$V^{*}(h)=\sqrt{ 2\gamma^{*}\Bar{g}_{sup}}\|K\|_{2}( \|K\|_1 + 1)( 1 + \delta^{*}_{n} ) \frac{\left(\log n\right)^{-\frac{1}{2}}}{\left(nh\right)^\frac{1}{2}}$,
where $\delta^*_n =(\log n)^{-\frac{1}{2}}$ and $\Bar{g}_{sup}$ denotes the maximum value of the histogram constructed from the sample $\{X_{i}\}_{i=1}^n$.

\begin{remark}
    In equation \eqref{CExpresionA*}, $\hat{g}_{h,h'}$ was used instead of $\hat{g}_{h\vee h'}$. Subsequently, an upper bound for $\hat{J}(h) = \frac{1}{n^2}\sum_{i=1}^{n} K^2_{h}(x-X_{i})$ was derived and substituted into equation \eqref{CExpresionV*}. Finally, this latter expression was further manipulated in order to obtain a form of $V^{*}(h)$ equivalent to that given in equation \eqref{CExpresionV}.
\end{remark}
The bandwidth $\hat{h}_j$ is selected from the collection of bandwidths $\mathcal{H}$ by taking
\begin{equation}\label{Chxi}
\hat{h}_{j}=\mbox{arg}\hspace{-0.05cm}\min_{h\in\mathcal{H}}\{A(h,x_{j})+V(h,x_{j})\}
\end{equation}
for each $j=1,\cdots,s$, where
$A(h,x_{j})= \max_{h'\in\mathcal{H}}\{|\hat{m}_{h,h'}(x_{j})-\hat{m}_{h'}(x_{j})|-V(h',x_{j})\}_+$,
with
$\hat{m}_{h,h'}(x_{j})=\frac{1}{n}\sum_{i=1}^nK_{h}\ast K_{h'}\left(x_{j}-X_{i}\right)Y_i$
and
$V(h,x_{j})=\sqrt{ 2\gamma A_{1}(j)}( \|K\|_1 + 1)( 1 + \delta_{n} ) \frac{\left(\log n\right)^\frac{1}{2}}{\left(nh\right)^\frac{1}{2}}$,
where
$\delta_n =(\log n)^{-\frac{1}{5}}$,
$A_{1}(j) = ((\hat{r}_{\sup}(j))^2+\hat{\sigma}^2)\hat{g}_{sup}\|K\|_2^2$,
$\hat{g}_{sup}=\max\{\hat{g}_{\hat{h}_j^*}(x_{j})| j=1,\cdots,s\}$,
$\hat{r}_{sup}(j)= \max\left\{|Y_{i}|\mathbf{1}_{]x_{j}-0.5,x_{j}+0.5[}(X_i)| i=1,\cdots,n \right\}$,
and
$\hat{\sigma}^2 =\frac{1}{n-1}\sum_{i=1}^n \left(Y_{i}-\tilde{r}_{\Bar{h}}(X_i)\right)^2$,
where $\Tilde{r}_{\Bar{h}}(\cdot)$ denotes the Nadaraya--Watson estimator associated with bandwidth $\Bar{h}$ obtained by the plug-in method proposed by Ruppert, Sheather, and Wand (1995).
\begin{figure}[htp!]
    \begin{center}
        \includegraphics[width=9.9cm]{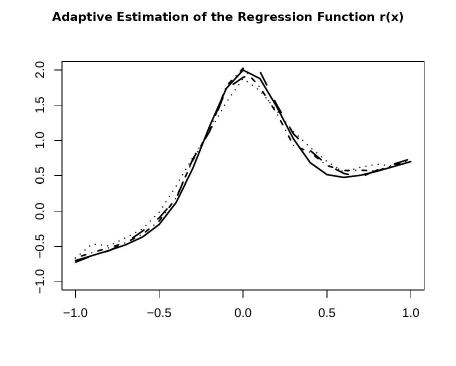}
    \end{center}
   \caption{Three estimators of the regression function $r$ are displayed for random samples $\{(X_{i},Y_{i})\}_{i=1}^{n+2q}$, with $\sigma=0.5$ and $q=100$. When $n=1000$, the estimator is represented by dots; when $n=2000$, by long dashes; and when $n=5000$, by dots and short dashes. The true regression function $r$ is also shown as a solid line.}
    \label{CGrafEstima}
\end{figure}

Figure~\ref{CGrafEstima} displays the regression function $r$ together with three estimators of the regression function obtained from random samples $\{(X_{i},Y_{i})\}_{i=1}^{n+2q}$, with $n=1000$, $2000$, and $5000$, taking $q=100$ and $\sigma=0.5$.\\

The parameters $\gamma^*$ and $\gamma$ are calibration parameters that control the penalty terms used by the GL bandwidth selection procedure to estimate the functions $g(\cdot)$ and $m(\cdot)$, respectively. Their values are determined through a data-driven calibration strategy described in Appendix~\ref{Calibra}. In particular, a grid of candidate values for $\gamma^*$ and $\gamma$ is considered, and the corresponding adaptive estimators are computed for each candidate value. The selection is performed using the calibration sample $\{(X_i,Y_i)\}_{i=n+q+1}^{n+2q}$, which is separated from the estimation sample by a time lag of length $q$ in order to attenuate the effects of serial dependence. Empirical criteria are constructed on this sample to select $\gamma^*$ and $\gamma$: for $\gamma^*$, the estimated density is compared with a reference histogram, while for $\gamma$, the relationship $m(x)\approx Y\,g(x)$ is used to assess the quality of the estimation of $m(\cdot)$. In this way, the smoothing level of the GL estimators is selected automatically from the data. The choice of the calibration grids considered in the numerical study and described in Appendix~\ref{Calibra} was motivated by the empirical observation of a transition between overfitting and underfitting regimes as $\gamma$ and $\gamma^{*}$ vary, a phenomenon reminiscent of the minimal-penalty results for Goldenshluger--Lepski procedures established by Lacour and Massart \cite{lacour2016minimalpen}. Although a theoretical analysis of this behavior is beyond the scope of the present work, its study could provide relevant insight into the theoretical behavior of the proposed adaptive estimator. Moreover, a rigorous characterization of these phenomena in the context of weakly dependent data could lead to a natural extension of the minimal-penalty theory developed for the classical setting of independent observations.

\subsection{Comparison of Global and Local Empirical Errors for Different Sample Sizes and Values of $\sigma$} \label{CEstErrores}

To assess the performance of the regression function estimator, $N=500$ replications of the sample $\{(X_{i},Y_{i})\}_{i=1}^{n+2q}$ are generated that allow to obtain $N$ calibrated estimators
$\left\{\hat{r}^{(k)}(x_{j})\right\}_{j=1}^s$, $k=1\ldots,N$.
 The mean squared error ($MSE$) and the mean integrated squared error ($MISE$) are estimated as follows:

$$\hat{MSE}(x_{j})=\frac{1}{N}\sum_{k=1}^{N}\left[ \left(r(x_{j})-\hat{r}^{(k)}(x_{j})\right)^2\right] \text{ for } j=1,\cdots,s$$

$$\hat{MISE}=\frac{1}{N}\sum_{k=1}^{N}I_k$$
with $I_{k}=\sum_{j=1}^{s-1}0.1\left(r(x_{j})-\hat{r}^{(k)}(x_{j})\right)^2$ for $k=1,\cdots,N$.\\

\begin{figure}[htp!] 
	\centering 
	\begin{subfigure}
	\centering
	\includegraphics[width=7.66cm]{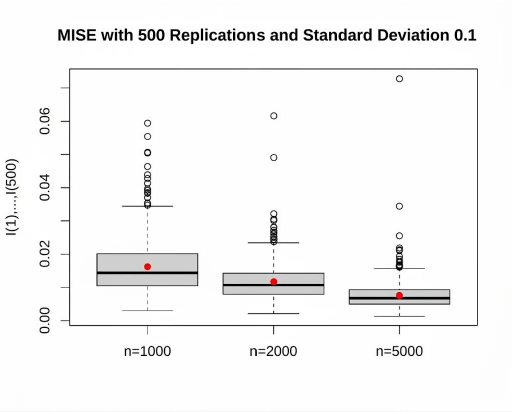}  
	\end{subfigure}
	\begin{subfigure}
	\centering
	\includegraphics[width=7.66cm]{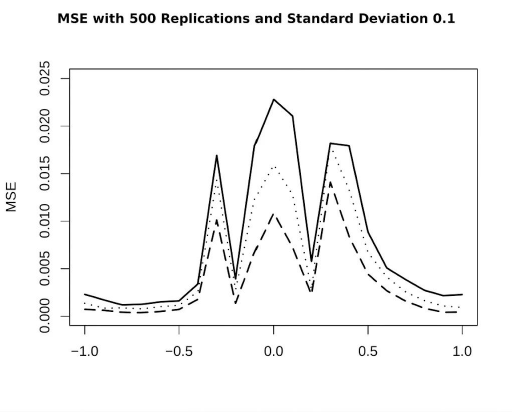}  
	\end{subfigure}
	\caption{The left panel displays three boxplots of $I_k$ for $k=1,\cdots,N$, corresponding to sample sizes $n=1000$, $2000$, and $5000$ from left to right, respectively. A point representing the estimated $\hat{MISE}$ is superimposed on each boxplot. The right panel shows three representations of $\hat{MSE}(x_{j})$ for $j=1,\cdots,s$: a solid line for $n=1000$, a dotted line for $n=2000$, and a dashed line for $n=5000$. In both panels, $\sigma=0.1$.}
    \label{CMSEyMISE01}
\end{figure}

In Figure~\ref{CMSEyMISE01}, three boxplots of $I_k$ for $k=1,\cdots,N$ are displayed, corresponding to sample sizes $n=1000$, $2000$, and $5000$ from left to right, respectively. In addition, a point representing the estimated $\hat{MISE}$ is superimposed on each boxplot (left panel). The right panel displays three representations of $\hat{MSE}(x_{j})$, for $j=1,\cdots,s$: a solid line for $n=1000$, a dotted line for $n=2000$, and a dashed line for $n=5000$. In both panels, $\sigma=0.1$. 

An examination of the $\hat{MSE}$ plot for $\sigma=0.1$ in Figure~\ref{CMSEyMISE01} (right panels) shows that the values of $\hat{MSE}$ are larger at points of the domain of $r$ close to its maximum and near the inflection region. This pattern is observed for all sample sizes considered, namely $n=1000$, $2000$, and $5000$. As expected, the values of $\hat{MSE}$ decrease as the sample size increases. Similar behavior has been observed for the $\hat{MSE}$ corresponding to $\sigma=0.5$ and $\sigma=1$. In these cases. as expected, the values of $\hat{MSE}$ increase as the value of $\sigma$ increases.

The boxplots of the values $I_{k}$ for $k=1,\cdots,N$ with $\sigma=0.1$ shown in Figure~\ref{CMSEyMISE01} (left panels) indicate that, as the sample size increases from $1000$ to $2000$ and then to $5000$, the dispersion of $\{I_k\}_{k=1}^N$ decreases. The number of outliers also decreases, and the value of $\hat{MISE}$ becomes progressively closer to the median of $\{I_k\}_{k=1}^N$, although the median remains generally smaller than $\hat{MISE}$. When $\sigma=0.5$ and $\sigma=1$, the dispersion of $\{I_k\}_{k=1}^N$ increases as $\sigma$ increases. 

Finally, Table~\ref{CT4tab:MISEHHg}, shows the $\hat{MISE}$ for different values of $n$ and $\sigma$. As expected, it decreases with increasing sample size $n$ and increases with increasing noise level $\sigma$.
\begin{table}[bt] 
\begin{center}
\begin{tabular}{lccc}
 & $\boldsymbol{n=1000}$ & $\boldsymbol{n=2000}$ & $\boldsymbol{n=5000}$\\
$\boldsymbol{\sigma=0.1}$ & 0.01624254 & 0.01177675 & 0.007650944\\
$\boldsymbol{\sigma=0.5}$ & 0.04625373 & 0.03409679 & 0.023852842\\
$\boldsymbol{\sigma=1}$ & 0.13286962055 & 0.10837311 & 0.073751940\\
\hline  
\end{tabular}
\caption{Estimated $\hat{MISE}$ values for different values of $\sigma$ and sample sizes in the case of unknown $g$.}
\label{CT4tab:MISEHHg}
\end{center}
\end{table}


 \section*{Conflict of interest}

The authors declare no potential conflict of interests.

\section*{Declaration of generative AI use}
AI has been used for Copyediting and language refinement.

\newpage
\appendix

\section{Technical Results} \label{CT4Demostraciones}
\vspace{0.5cm}

\subsection{Constants} \label{CT4Constantes}
We recall some constant values given previously, and below we declare the constants that will be used in the development of the proofs.

\begin{align*}
& A_0 =\frac{L_1}{\mathit{l}_1!}\int_{\mathbb{R}}|u|^{\beta_1}|K(u)|du, \quad A_{1}=((r_{\sup})^2 +\sigma^2 )g_{\sup}\|K\|_2^{2}=B_{1}+B_{2},\\
&
A_{2} = \frac{2A_3}{|\log a|} + 20A_4, \quad A_3 = (r_{\sup})^2 \|K\|_1^2 \left(Q+(g_{\sup})^{2}\right), \quad A_4 =\frac{4}{1-a}(r_{\sup})^2 \|K\|_{\infty}^2 ,\\
&
C=A_{0}^{2}+A_{1}+A_{2}, \quad \mathfrak{A}_0 =\frac{L_2}{\mathit{l}_2 !}\int_{\mathbb{R}}|u|^{\beta_2}|K(u)|du, \quad \mathfrak{A}_{1}=g_{\sup}\|K\|_2^{2},\\
&
\mathfrak{A}_{2} = \frac{2\mathfrak{A}_3}{|\log a|} + 20\mathfrak{A}_{4}, \quad \mathfrak{A}_3 = \|K\|_1^2 \left(Q+(g_{\sup})^{2}\right), \quad\mathfrak{A}_4 =\frac{4}{1-a}\|K\|_{\infty}^{2},\\
& \mathfrak{C}=\mathfrak{A}_{0}^{2}+\mathfrak{A}_{1}+\mathfrak{A}_{2}, \quad B_{1} =\sigma^2 g_{\sup}\|K\|_2^{2}, \quad B_{2} =(r_{\sup})^2 g_{\sup} \|K\|_2^{2},\\
&
B_{3} =\left( r_{\sup} + \frac{2\sigma}{\sqrt{2\pi}}\right)g_{\sup}\|K\|_{1}, \quad B_{4} = \left( r_{\sup} + \frac{2\sigma}{\sqrt{2\pi}}\right)^2 Q \|K\|_1^{2},\\
&
B_{5} = \left( r_{\sup} + \frac{2\sigma}{\sqrt{2\pi}}\right)^2 \|K\|_1^2 \left(Q+(g_{\sup})^2 \right),  \quad  B_{6} = 2\|K\|_{\infty}, \quad B_{7} = (4B_{5} B_{6} (\sigma+r_{\sup}))^{\frac{2}{3}},\\
&
B_{8} = \frac{2B_{5}}{|\log a|} +\frac{18B_{7}}{1-a^\frac{1}{3}}, \quad  B=\frac{12(B_{6}\vee \sqrt{B_7})}{1 \wedge |\log a|}\left( \frac{2^{7}3B_{7}a^{-\frac{1}{3}}}{A_{1}(1 \wedge |\log a |)} \vee 1\right),\\
&
{B_{9}= \frac{2\|K\|^2_{\infty}}{(2\pi)^{\frac{1}{4}}}\left(2(r_{\sup})^4 +3\sigma^4\right)^{\frac{1}{2}}}, \quad B_{10} = \frac{2^4 3}{1-e^{-\frac{\gamma-2}{2}}}(A_{1} + B_{8}) + \frac{5!2^{18} 3 B^2}{1-e^{-\frac{\gamma-2}{2}}} (\sigma + r_{\sup})^{2},\\
&
A_5=B_9+2B_{10}, \quad A_{6}=1+2\|K\|_{1}, \quad A_{7}=\sqrt{A_{1}+A_{2}}+2\sqrt{2\gamma A_1}(1+\|K\|_{1}),\\
&
A_{8}=2\sqrt{A_5}(1+\|K\|_{1}), \quad C^{*}=A_{0}A_{6}+2e^{\frac{1}{2}}A_{7}+A_{8},\\
&
\Bar{F}_{1}=  2B_{11}(r_{sup}g_{sup})^2 + 2\|K\|^2_{\infty}B_{11}^{1/2}(2(2r_{sup}^4+3\sigma^4))^{1/2}, \quad \Bar{F}_{2}=2B_{11}\left\{\left(\|K\|^2_{\infty} \vee 1\right) + (g_{\sup})^{2}\right\},\\
&
F_{1}=\frac{8}{(g_{inf})^2}, \quad F_{2}=\frac{8(r_{sup})^2}{(g_{inf})^2}, \quad F_{3}=2\Bar{F}_{1} + 2(r_{sup})^2\Bar{F}_{2}, \quad \mathfrak{D}=F_{1}C+F_{2}\mathfrak{C}+F_3,\\
&
\mathfrak{B}=\left(F_{1}({C^*})^2 +F_{2}(C^{**})^{2} +F_{3}\right)^{1/2}.
\end{align*}

\subsection{Statements of technical results} \label{CT4ERT}

\begin{lemma} \label{CT4L1}
Under hypotheses $(H_1)$, $(H_2)$, $(H_4)$, and $(H_5)$, we have that {for $0\le h\le 2h_{max}$}:
\begin{itemize}
    \item[i)] $\mathbb{E}\left[\left( K_h (x-X_i)\varepsilon_i \right)^2\right] \leq \frac{B_{1}}{h}$ 
    \item[ii)] $\mathbb{E}\left[ \left( K_h (x-X_i)r(X_i) \right)^2 \right] \leq \frac{B_{2}}{h}$ 
    \item[iii)] $\mathbb{E}\left[ |K_h(x-X_i)K_h(x-X_j)| \right] \leq Q\|K\|_1^2$
    \item[iv)] $\mathbb{E}\left[ |K_h(x-X_i)r(X_i)| \right] \leq r_{\sup} g_{\sup} \|K\|_1$,
\end{itemize}
where $B_{1} =\sigma^2 g_{\sup}\|K\|_2^2$ and $B_{2} =(r_{\sup})^2 g_{\sup} \|K\|_2^2$.

\end{lemma}

\begin{lemma} \label{CT4L2}
 For each{ $0\le h\le h_{max}$}, $a\in]0,1[$, and $n\geq 4$, it holds that,
 $$\frac{1}{h}a^{\left\{\frac{1}{h|\log a|(\log n)^\frac{1}{2}}\right\}} \leq \frac{10}{(\log n)^\frac{5}{2}}$$
\end{lemma}
\begin{proposition} \label{CT4P1}
 We assume that the $X_i$ satisfy the hypothesis $(H_3)$. Let $k,u,v\in\mathbb{N}$ such that ${i_{1:u}}\in\mathbb{Z}^u$, ${j_{1:v}}\in\mathbb{Z}^v$, and $i_1 \leq \cdots \leq i_u < i_u +k \leq j_1 \leq \cdots \leq j_v$.
  If $G_u :\mathbb{R}^u \times\mathbb{R}^u \rightarrow \mathbb{R}$ and $G_v :\mathbb{R}^v \times\mathbb{R}^v \rightarrow \mathbb{R}$ are functions such that $\|G_u\|_{\infty}<\infty$ and $\|G_v\|_{\infty}<\infty$, then,
 $$|cov\left(G_u(X_{i_{1:u}},\varepsilon_{i_{1:u}} ),G_v (X_{j_{1:v}} ,\varepsilon_{j_{1:v}} )\right)|\leq \Psi (u,v,G_u ,G_v ) a^k$$
 where $\Psi(u,v,G_u , G_v )=4\|G_u \|_{\infty}\|G_v\|_{\infty}$, with $\|G_u\|_{\infty}=\sup_{(x,y)\in\mathbb{R}^u\times\mathbb{R}^u}|G_u (x,y)|$ and $\|G_v\|_{\infty}=\sup_{(x,y)\in\mathbb{R}^v\times\mathbb{R}^v}|G_v (x,y)|$.
 \end{proposition}
{Lemma~\ref{CT4L2} and Proposition~\ref{CT4P1} are proved in \citet{mpadrino2024adaptiregre}.} As in \citet{mpadrino2024adaptiregre}, next we will work with specific functions $G_u$. More precisely, for $u\in\mathbb{N}$ and ${i_{1:u}}\in\mathbb{Z}^u$, we consider $G_u :\mathbb{R}^u \times\mathbb{R}^u \rightarrow \mathbb{R}$ given by
 $$G_u (X_{i_{1:u}} ,\varepsilon_{i_{1:u}} )=\prod_{k=1}^u\left\{G(X_{i_k},\varepsilon_{i_k})-\mathbb{E}\left(G(X_{i_k},\varepsilon_{i_k})\right)\right\}$$
 where for $i\in\mathbb{N}$,
\begin{equation}\label{newG}
G(X_{i},\varepsilon_{i}) =\frac{1}{n}(r(X_{i})+\varepsilon_{i})K_h (x-X_i)\mathbf{1}_{\{|(r(X_{i})+\varepsilon_{i}|\leq M_n \}},
\end{equation}
where $M_n = \sigma \log n +r_{sup}$. The variables $G(X_{i},\varepsilon_{i})$ satisfy the following results.
  
  \begin{lemma} \label{CT4L3}
  Under hypotheses $(H_1)$, $(H_2)$, $(H_4)$, and $(H_5)$, {for $0\le h\le 2h_{max}$}, we have that: $\forall i,k \in \mathbb{N}$
 \begin{enumerate}
     \item[i)] $|\mathbb{E}\left[G(X_i ,\varepsilon_i ) \right]| \leq \frac{B_{3}}{n}$,
     \item[ii)] $|\mathbb{E}\left[ G(X_i ,\varepsilon_i )G(X_{i+k} ,\varepsilon_{i+k}) \right]| \leq \frac{B_{4}}{n^2}$,
     \item[iii)] $|Cov\left( G(X_i ,\varepsilon_i ), G(X_{i+k} ,\varepsilon_{i+k}) \right)| \leq \frac{B_{5}}{n^2}$,
     \item[iv)] $|G(X_{i},\varepsilon_{i})| \leq \frac{B_{6}}{2} \frac{M_n}{nh}$
 \end{enumerate}
 where $$B_{3} =\left( r_{\sup} + \frac{2\sigma}{\sqrt{2\pi}}\right)g_{\sup}\|K\|_1, \quad B_{4} = \left( r_{\sup} + \frac{2\sigma}{\sqrt{2\pi}}\right)^2 Q \|K\|_1^2 $$ 
 $$B_{5} = \left( r_{\sup} + \frac{2\sigma}{\sqrt{2\pi}}\right)^2 \|K\|_1^2 \left(Q+(g_{\sup})^2 \right) \quad \text{and} \quad B_{6} = 2\|K\|_{\infty}.$$
 \end{lemma}

 \begin{proposition} \label{CT4P2}
Let $ u,v,k \in \mathbb{N}$, and let $(i_1 , \cdots ,i_u ,j_1 , \cdots ,j_v)\in \mathbb{Z}^{u+v}$ be such that $i_1 \leq \cdots \leq i_u < i_u +k \leq j_1 \leq \cdots \leq j_v$. Under hypotheses $(H_1)$, $(H_2)$, $(H_3)$, $(H_4)$, and $(H_5)$, we have $\forall h \in \mathcal{H}$
 $$\gamma_h (u,v)=|cov\left(G_u (X_{i_{1:u}} ,\varepsilon_{i_{1:u}}), G_v (X_{j_{1:v}} ,\varepsilon_{j_{1:v}})\right)| \leq \phi (u,v)(D_1(h))^{u+v-2}D_2 (h)  a^{\frac{k}{3}}$$
 where $\phi(u,v)= u+v+uv$,
 $$D_1 (h)=\frac{B_{6} M_n}{nh}, \quad D_2 (h)=\frac{B_{7}}{n^2 h}\quad \text{and \ } B_{7} = (4B_{5} B_{6} (\sigma+r_{\sup}))^{\frac{2}{3}}.$$
 \end{proposition}

 \begin{proposition} \label{CT4TDBernstein}
\textbf{Bernstein's Inequality}. Let the process be $\mathbb{Z}=(Z_i)_{i\in\mathbb{N}}$, where $Z_i =G(X_i,\varepsilon_i)-\mathbb{E}[G(X_i,\varepsilon_i)]$ and $S_n=\sum_{i=1}^n Z_i$.
Under hypotheses $(H_1)$, $(H_2)$, $(H_3)$, $(H_4)$, and $(H_5)$, for $h\in\mathcal{H}$ and $\forall t \geq 0$, we have
\begin{equation} \label{CT4DBernstein}
    \mathbb{P}(|S_n | \geq t) \leq \exp \left(-\frac{\frac{t^2}{2}}{\mathcal{A}_n +\mathcal{B}_n^\frac{1}{3}t^\frac{5}{3}}\right)
\end{equation}
where $\mathcal{A}_n=\frac{A_{1}}{nh} + B_{8} \frac{(\log n)^{-\frac{1}{2}}}{nh}$ is an upper bound of $Var(S_{n})$, and $\mathcal{B}_n =B\frac{M_n}{nh}$, with
$$B_{8} = \frac{2B_{5}}{|\log a|} +\frac{18B_{7}}{1-a^\frac{1}{3}},\quad  B=\frac{12(B_{6}\vee \sqrt{B_7})}{1 \wedge |\log a|}\left( \frac{2^{7}3B_{7}a^{-\frac{1}{3}}}{A_{1}(1 \wedge |\log a |)} \vee 1\right).$$
\end{proposition}

 \begin{corollary} \label{CT4CDBernstein} Under the hypotheses of Proposition \ref{CT4TDBernstein}, if $\lambda(u)=(2\mathcal{A}_n u)^\frac{1}{2}+\mathcal{B}_n (2u)^3$ for each $u\geq 0$ and we perform the change of variable $t=\lambda(u)$ in Bernstein's inequality \eqref{CT4DBernstein}, then the inequality can be rewritten as follows:
 
 $$ \mathbb{P}(|S_n | \geq \lambda(u)) \leq \exp \left(-\frac{u}{2}\right)$$
 for each $u \geq 0$.
 \end{corollary}
{Proofs of Propositions~\ref{CT4P2} and \ref{CT4TDBernstein}, and Corollary~\ref{CT4CDBernstein} follow exactly the proof of Proposition 6 and 7, and Corollary 1 respectively, of \citet{mpadrino2024adaptiregre} using the function $G$ defined in \eqref{newG}.}

\begin{lemma} \label{CT4MomenY^2,4,8}
For $i\in\{1,\ldots,n\}$ and {$0<h\leq 2h_{\max}$}, we have:
{\begin{itemize}
     \item[i)] $\E[Y_i^2\mathbf{1}_{\{|X_i-x|\leq h\}}]\leq (r_{\sup})^2 +\sigma^2$
     \item[ii)]  $\E[Y_i^4\mathbf{1}_{\{|X_i-x|\leq h\}}]\leq 2\left(2(r_{\sup})^4 +3\sigma^4\right)$
 \end{itemize}}
\end{lemma}

\begin{lemma} \label{CT4P|Y|>Mn}
For $i\in\{1,\ldots,n\}$ and {$0<h\leq 2h_{\max}$}, we have: $$\mathbb P[(|Y_i|>M_n)\cap \{|X_i-x|\leq h\}]\leq \frac{2}{\sqrt{2\pi}}\frac{e^{-\frac{\log^2 n}{2}}}{\log n}.$$
\end{lemma}

\begin{lemma} \label{CT4MomenEta1,2y4}
    Let $\eta_i = Y_i K_h (x-X_i)\mathbf{1}_{\{|Y_i|>M_n\}}$. 
    Then, for {$0<h\leq 2h_{\max}$}
    \begin{itemize} 
        \item[i)] {$\E[|\eta_i|]\leq\frac{\|K\|_{\infty}\sqrt{2}}{(2\pi)^{\frac{1}{4}}}((r_{\sup})^{2}+\sigma^2 )^{\frac{1}{2}}\frac{e^{-\frac{(\log n)^2}{4}}}{h(\log n)^{\frac{1}{2}}}$}
        \item[ii)] $\E[\eta_i^2]\leq B_{9}\frac{e^{-\frac{(\log n)^2}{4}}}{h^{2}(\log n)^{\frac{1}{2}}}$
        where {$B_{9}=\frac{2\|K\|^2_{\infty}}{(2\pi)^{\frac{1}{4}}}\left(2(r_{\sup})^4 +3\sigma^4\right)^{\frac{1}{2}}$.}
    \end{itemize}
\end{lemma}

\begin{proposition} \label{CT4ProOmega_h^c}
    For  {$h\in[h_{min}^*,h_{max}^*]$}, we consider the event $\Omega_h =\{\hat{g}_{h}(x)\geq \frac{g_{\inf}}{2}\}$. 
    Assume that $g\in\Sigma(\beta_2 ,L_2 )$ with $\beta_2 >0$, $L_2>0$, and $K$ is a kernel of order $\ord \geq l_2 =[\beta_2 ]$ with $\int_{\mathbb{R}}|u|^{\beta_2}|K(u)|du<\infty$. Under $(H_1)$, $(H_2)$, $(H_3)$, and $(H_5)$, we have
        $$\mathbb{P}(\Omega_h^c )\leq \frac{B_{11}}{n^{13}}$$ 
    where $B_{11}$ is a positive constant that depends on the kernel $K$, $|\log(a)|$, $\{\alpha_{k}(\mathbb{X}) \}_{k\geq 0}$, $g_{\sup}$, and $Q$.
\end{proposition}

\begin{proposition} \label{CT4ProOmega_Hath^c}
    Under the hypotheses of Proposition \ref{CT4ProOmega_h^c}. If $\hat{h}^* =argmin_{h\in\mathcal{H}^*} \{A^* (h,x)+V^* (h)\}$, then

    $$\mathbb{P}(\Omega_{\hat{h}^*}^c )\leq \frac{B_{11}}{n^{12}}.$$   
\end{proposition}

\subsection{Proofs of technical results} 

\subsubsection{Proof of Lemma \ref{CT4L1}}
i) We have
$$
    \mathbb{E}\left[\left(K_h (x-X_i)\varepsilon_i \right)^2\right] =
    \frac{\sigma^2}{h^2} \mathbb{E}\left[K^2 \left(\frac{x-X_i}{h}\right) \right]=
    \frac{\sigma^2}{h^2} \int_{\mathbb{R}}K^2 \left(\frac{x-t}{h}\right) g(t)dt\le 
    \frac{B_{1}}{h}.$$
ii) We have
\begin{eqnarray}
   \nonumber
   \mathbb{E}\left[ \left( K_h (x-X_i)r(X_i) \right)^2 \right] &=&  \mathbb{E} \left[\frac{1}{h^2} K^2 \left(\frac{x-X}{h}\right)r^{2}(X)\right]\\
   \nonumber
   &=&
   \frac{1}{h^2} \int_{\mathbb{R}} K^2 \left(\frac{x-t}{h}\right)r^2(t )g(t)dt\\
   \nonumber
   &\leq&
   \frac{(r_{\sup})^2 g_{\sup}}{h^2} \int_{\mathbb{R}} K^2 (u)hdu=\frac{B_{2}}{h} .
\end{eqnarray}
iii) {See Lemma 2 in \citet{mpadrino2024adaptiregre}}\\
iv) We have
\begin{eqnarray*}
   \mathbb{E}\left[ |K_h(x-X_i)r(X_i)| \right]  &\leq& r_{\sup} \mathbb{E}\left[ |K_h(x-X_i)| \right]\\
   &=&
    r_{\sup} \int | K_h(x-t)|g(t)dt\\
   &\leq&
 r_{\sup} g_{\sup}\int \frac{1}{h}\left| K\left(\frac{x-t}{h}\right) \right|dt= r_{\sup}g_{\sup}\|K\|_1.
\end{eqnarray*}

\subsubsection{Proof of Lemma \ref{CT4L3}}

i) We have
 \begin{eqnarray*}
 |\mathbb{E}(G(X_i,\varepsilon_i))| &\leq& \mathbb{E}\left\{\left|\frac{1}{n}(r(X_i)+\varepsilon_i)K_h (x-X_i)\mathbf{1}_{\{|r(X_i)+\varepsilon_i|\leq M_n \}}\right|\right\}\\
  &\leq& 
  \frac{1}{n}\mathbb{E}\left\{|r(X_i)+\varepsilon_i||K_h (x-X_i)|\right\}\\
  &\leq& 
  \frac{1}{n}\mathbb{E}\left\{|r(X_i)K_h (x-X_i)|\right\}+\frac{1}{n}\mathbb{E}\left\{|\varepsilon_i K_h (x-X_i)|\right\}\\
  &\leq& 
  \frac{1}{n}r_{\sup}g_{\sup}\|K \|_1 +\frac{1}{n}\frac{2\sigma}{\sqrt{2\pi}}g_{\sup}\|K \|_1=\frac{B_{3}}{n}
 \end{eqnarray*}
by inequality $iv)$ of Lemma \ref{CT4L1}.\\

ii) Denoting $A(X_i)=\frac{1}{n}K_h (x-X_i)$, by using inequality (iii) of Lemma \ref{CT4L1}, we have

\begin{eqnarray*}    \lefteqn{|\mathbb{E}\left[G(X_i,\varepsilon_i)G(X_{i+k},\varepsilon_{i+k})\right]|}\\
    &\leq& 
    \mathbb{E}\left[|r(X_i)+\varepsilon_i||A(X_i)||r(X_{i+k})+\varepsilon_{i+k}||A(X_{i+k})|\right]\\
     &\leq& 
    \mathbb{E}\left[|r(X_i)r(X_{i+k})A(X_i)A(X_{i+k})|\right] + \mathbb{E}\left[|r(X_i)\varepsilon_{i+k}A(X_i)A(X_{i+k})|\right]\\
    &+&
    \mathbb{E}\left[|\varepsilon_i r(X_{i+k})A(X_i)A(X_{i+k})|\right] + \mathbb{E}\left[|\varepsilon_i\varepsilon_{i+k}A(X_i)A(X_{i+k})|\right]\\
     &\leq& 
    (r_{\sup})^2\mathbb{E}\left[|A(X_i)A(X_{i+k})|\right] + r_{\sup}\mathbb{E}[|\varepsilon_{i+k}|]\mathbb{E}\left[|A(X_i)A(X_{i+k})|\right]\\
    &+&
    r_{\sup}\mathbb{E}[|\varepsilon_i|]\mathbb{E}\left[|A(X_i)A(X_{i+k})|\right] + \mathbb{E}[|\varepsilon_i|]\mathbb{E}[|\varepsilon_{i+k}|]\mathbb{E}\left[|A(X_i)A(X_{i+k})|\right]\\
     &=& 
    \left\{(r_{sup})^2 + 2r_{\sup}\frac{2\sigma}{\sqrt{2\pi}}+\left(\frac{2\sigma}{\sqrt{2\pi}}\right)^2\right\}\mathbb{E}\left[|A(X_i)A(X_{i+k})|\right]\\
    &\leq&
     \left(r_{\sup} + \frac{2\sigma}{\sqrt{2\pi}}\right)^2 \frac{1}{n^2}\mathbb{E}\left[|K_h(x-X_i)K_h(x-X_{i+k})|\right]\\
      &\leq&
     \left(r_{\sup} + \frac{2\sigma}{\sqrt{2\pi}}\right)^2 Q\|K\|_1^2\frac{1}{n^2}\\
     &=&
     \frac{B_{4}}{n^2}.
\end{eqnarray*}

iii) Now using (i) and (ii), we have
\begin{eqnarray*}
\lefteqn{|Cov\left( G(X_i ,\varepsilon_i ), G(X_{i+k} ,\varepsilon_{i+k}) \right)|}\\
&\leq&
    |\mathbb{E}[G(X_i,\varepsilon_i)G(X_{i+k},\varepsilon_{i+k})]|+|\mathbb{E}[G(X_i,\varepsilon_i)]||\mathbb{E}[G(X_{i+k},\varepsilon_{i+k})]|\\
    &\leq&
    \frac{B_{4}}{n^2}+\frac{B_{3}}{n}\frac{B_{3}}{n}=
    \frac{B_{5}}{n^2}.
 \end{eqnarray*}
Inequality (iv) follows immediately using hypothesis $(H_5)$.

\subsubsection{Proof of Lemma \ref{CT4MomenY^2,4,8}}

{
i) We have using the independence of the $X_i$ and the $\varepsilon_i$ and that $\E(\varepsilon_i)=0$
\begin{align*}\E[Y_i^2\mathbf{1}_{\{|X_i-x|\leq h\}}]&=
     \E[(r^2(X_i)+\varepsilon_i^2)\mathbf{1}_{\{|X_i-x|\leq h\}}]\\
     &\leq (r_{\sup})^2 +\sigma^2
     \end{align*}
     where the last line is a consequence of $h\le 2h_{max}$.\\
     ii) Using the same arguments, we deduce that
     \begin{align*}\E[Y_i^4\mathbf{1}_{\{|X_i-x|\leq h\}}]&=
     \E[(r^4(X_i)+6r^2(X_i)\varepsilon_i^2+\varepsilon_i^4)\mathbf{1}_{\{|X_i-x|\leq h\}}]\\
     &\leq \left((r_{\sup})^4 +6(r_{\sup})^2\sigma^2+3\sigma^4\right)\\
     &\leq 2\left(2(r_{\sup})^4 +3\sigma^4\right)
     \end{align*}
}

\subsubsection{Proof of Lemma \ref{CT4P|Y|>Mn}}
We have using that {$h\le 2h_{max}$}
{\begin{eqnarray*}
\mathbb{P}\left( \{|Y_i |> M_n\}\cap \{|X_i-x|\leq h\}\right) &=& \mathbb{P}\left( \{|r(Xi)+\varepsilon_i |> M_n\}\cap \{|X_i-x|\leq h\}\right)\\
&\leq&
\mathbb{P}\left( \{r_{\sup}+|\varepsilon_i |> M_n\}\right)\\
&\leq&
\mathbb{P}\left( \left\{ |\varepsilon_i |> M_n - r_{\sup}\right\} \right)\\
&\leq&
2\mathbb{P}\left( \left\{\frac{\varepsilon_i}{\sigma}  > \frac{ M_n - r_{\sup}}{\sigma}\right\}\right)\\
&\leq&
2\frac{\phi\left( \log n\right)}{\log n}
\end{eqnarray*}}
substituting $M_n$ by its value and where $\phi(x)=\frac{1}{\sqrt{2\pi}}e^{-\frac{x^2}{2}}$. From the above, we have that
\begin{equation*}
    \mathbb{P}\left( \{|Y_i |> M_n\}\cap \{|X_i-x|\leq h\}\right) \leq \frac{2}{\sqrt{2\pi}}\frac{e^{-\frac{(\log n)^2}{2}}}{\log n}.
\end{equation*}

\subsubsection{Proof of Lemma \ref{CT4MomenEta1,2y4}}
i)
{\begin{eqnarray*}
    \E[|\eta_i|] &=& \E\left[ |Y_i\frac{1}{h}K\left(\frac{x-X_i}{h}\right)\mathbf{1}_{\{|Y_i|>M_n\}}|\right]\\
    &\leq&
    \frac{\|K\|_{\infty}}{h}\E\left[|Y_i|\mathbf{1}_{\{|Y_i|>M_n \}} \mathbf{1}_{\{|X_i-x|\leq h\}}\right]\\
    &\leq&
    \frac{\|K\|_{\infty}}{h}\left(\E[|Y_i|^{2}\mathbf{1}_{\{|X_i-x|\leq h\}}]\right)^{\frac{1}{2}}\left(\mathbb{P} \left[\{|Y_i|>M_n \}\cap \{|X_i-x|\leq h\}\right]\right)^{\frac{1}{2}}\\
    &\leq&
    \frac{\|K\|_{\infty}}{h}\left((r_{\sup})^{2}+\sigma^{2}\right)^{\frac{1}{2}}\left(\frac{2}{\sqrt{2\pi}}\frac{e^{-\frac{(\log n)^2}{2}}}{\log n} \right)^{\frac{1}{2}}\\
    &=&
    \frac{\|K\|_{\infty}\sqrt{2}}{(2\pi)^{\frac{1}{4}}}\left((r_{\sup})^{2}+\sigma^2\right)^{\frac{1}{2}}\frac{e^{-\frac{(\log n)^2}{4}}}{h(\log n)^{\frac{1}{2}}}
\end{eqnarray*}}
by part $(i)$ of Lemma \ref{CT4MomenY^2,4,8} and Lemma \ref{CT4P|Y|>Mn}.

ii) Using Cauchy-Schwarz again, part $(ii)$ of Lemma \ref{CT4MomenY^2,4,8}, and Lemma \ref{CT4P|Y|>Mn}, we obtain:
{\begin{eqnarray*}
    \E[\eta_i^2] &\le &
    \frac{\|K\|^2_{\infty}}{h^2}\left(\E[(Y_{i})^{4}\mathbf{1}_{\{|X_i-x|\leq h\}}]\right)^{\frac{1}{2}}\left(\mathbb{P} \left[\{|Y_i|>M_n \}\cap \{|X_i-x|\leq h\}\right]\right)^{\frac{1}{2}}\\
    &\leq&
    \frac{2\|K\|^2_{\infty}}{(2\pi)^{\frac{1}{4}}}\left(2(r_{\sup})^{4}+3\sigma^4\right)^{\frac{1}{2}}\frac{e^{-\frac{(\log n)^2}{4}}}{h^{2}(\log n)^{\frac{1}{2}}}
\end{eqnarray*}}

\subsubsection{Proof of Proposition \ref{CT4ProOmega_h^c}}
Let {$h\in[h_{min}^*,h_{max}^*]$. }We have
\begin{eqnarray} \label{CT4ProOmega^c11}
    \nonumber
    \mathbb{P}[\Omega_h^c] &=& \mathbb{P}\left[\hat{g}_{h}(x)<\frac{g_{\inf}}{2}\right]\\
    \nonumber
    &=&
    \mathbb{P}\left[ \hat{g}_{h}(x)-\E[\hat{g}_{h}(x)]+\E[\hat{g}_{h}(x)]-g(x)+g(x)<\frac{g_{\inf}}{2} \right]\\
    \nonumber
    &=&
    \mathbb{P}\left[\frac{1}{nh}\sum_{i=1}^n \left\{K\left(\frac{x-X_i}{h}\right)-\E\left[K\left(\frac{x-X_i}{h} \right)\right] \right\}+b(x)+g(x) < \frac{g_{\inf}}{2}\right]\\
    &=&
    \mathbb{P}\left[S_n < \frac{g_{\inf}}{2}-g(x)-b(x) \right]
\end{eqnarray}
where $S_{n}=\frac{1}{nh}\sum_{i=1}^n \left\{K\left(\frac{x-X_i}{h}\right)-\E\left[K\left(\frac{x-X_i}{h} \right)\right] \right\}$ and $b(x)=\E[\hat{g}_{h}(x)]-g(x)$.\\
Using \eqref{CT4pro:sesgoHatgh} of Proposition \ref{resume}, we have
$$|b(x)|\leq h^{\beta_{2}}\mathfrak{A}_{0}$$
and therefore,
\begin{equation*}
    -h^{\beta_{2}}\mathfrak{A}_0<-b(x)<h^{\beta_{2}}\mathfrak{A}_0
\end{equation*}
By the previous inequality and since $\frac{g_{\inf}}{2}-g(x)<-\frac{g_{\inf}}{2}$, we have
\begin{equation*} 
    \frac{g_{\inf}}{2}-g(x)-b(x)<-\frac{g_{\inf}}{2}+h^{\beta_{2}}\mathfrak{A}_0<-\frac{g_{\inf}}{2}+(h_{\max}^*)^{\beta_{2}}\mathfrak{A}_0
\end{equation*}

Thus, for $n$ sufficiently large, we have
\begin{equation*} 
    \frac{g_{\inf}}{2}-g(x)-b(x)<-\frac{g_{\inf}}{4}
\end{equation*}

By inequality \eqref{CT4ProOmega^c11}, it follows that
\begin{equation} \label{CT4ProOmega^c12}
    \mathbb{P}[\Omega_h^c] \leq \mathbb{P}\left[S_n <-\frac{g_{\inf}}{4} \right]\leq 
    \mathbb{P}\left[|S_n| >\frac{g_{\inf}}{4} \right]
\end{equation}
In Lemma 3 of \citet{bertin2017pointwise}, it is established that
\begin{equation} \label{CT4ProOmega^c13}
    \mathbb{P}\left[|S_{n}|\geq \lambda^{*}(t)\right] \leq e^{-\frac{t}{2}}
\end{equation}
for $t\geq 0$, where $\lambda^{*}(t)=\sigma_n^* (h)\sqrt{2t}+B_n^{*}(h)(2t)^3$, $\sigma_n^{*}(h) =\sqrt{ J_{n}+C_{6}^{*}\frac{(\delta_n^*)^2}{nh}}$, $\delta_n^{*}=(\log n)^{-\frac{1}{2}}$, $B_n^{*}(h)=\frac{C_7^*}{nh}$, $J_{n}(h)=\frac{1}{n}\int K_h^{2}(x_{0}-t)g(t)dt$, $C_6^*$ is a positive constant depending on the kernel $K$, $|\log(a)|$, $\{\alpha_{k}(\mathbf{X})\}_{k\geq 0}$, $g_{sup}$, and $Q$, and $C_7^*$ is a positive constant depending on the kernel $K$, $|\log(a)|$, $g_{sup}$, and $Q$.\\

By definition of $J_{n}(h)$, we have
\begin{equation*}
    J_{n}(h)\leq \frac{g_{\sup}}{nh^2}\int K^2\left(\frac{x_{0}-t}{h} \right)dt \leq \frac{g_{\sup}\|K\|_2^2}{nh}
\end{equation*}
By the bound of $J_{n}(h)$, the definitions of $B_n^{*}(h)$, $\delta_n^{*}$, $\sigma_n^{*}(h)$, and $\lambda^{*}(t)$, taking $t_{0}=26\log n$, we have
\begin{eqnarray*}
    \lambda^{*}(t_{0}) &\leq& \sqrt{\frac{g_{\sup}\|K\|_2^2}{nh}+C_{6}^*\frac{(\log n)^{-1}}{nh}}\sqrt{2t_0}+\frac{C_7^*}{nh}(2t_{0})^3\\
    &\leq&
    \sqrt{\frac{g_{\sup}\|K\|_2^2}{nh_{\min}^*}+C_{6}^*\frac{(\log n)^{-1}}{nh_{\min}^*}}\sqrt{2(26\log n)}+\frac{C_7^*}{nh_{\min}^*}(2(26\log n))^3\\
     &=&
     \sqrt{\frac{g_{\sup}\|K\|_2^2}{e^{(\sqrt{\log n}-1)}}+C_{6}^*\frac{(\log n)^{-1}}{e^{(\sqrt{\log n}-1)}}}\sqrt{52\log n}+\frac{C_7^*}{e^{(\sqrt{\log n}-1)}}(52\log n)^3
\end{eqnarray*}
It is clear that $\lambda^{*}(t_{0})\searrow0$ as $n\rightarrow\infty$ with $t_{0}=26\log n$; therefore, for $n$ sufficiently large, we have $\lambda^{\ast}(t_{0})<\frac{g_{\inf}}{4}$. This allows us to establish the following event inclusion relation:
$$\left[|S_n|>\frac{g_{\inf}}{4}\right]\subseteq\left[|S_n|>\lambda^{*}(t_{0})\right]$$
By the above and inequality \eqref{CT4ProOmega^c13}, we have
$$\mathbb{P}\left[|S_n|>\frac{g_{\inf}}{4} \right]\leq e^{-\frac{26\log n}{2}}=n^{-13}$$
for $n$ sufficiently large. By this last inequality and inequality \eqref{CT4ProOmega^c12}, we conclude that

$$\mathbb{P}[\Omega_h^c]\leq B_{11}n^{-13},$$
where $B_{11}$ depends on the kernel $K$, $|\log(a)|$, $\{\alpha_{k}(\mathbb{X}) \}_{k\geq 0}$, $g_{\sup}$, and $Q$.

\subsubsection{Proof of Proposition \ref{CT4ProOmega_Hath^c}}
We have
\begin{equation} \label{CT4ProOmega_Hath^c12}
    \nonumber   \mathbb{P}\left[\Omega^c_{\hat{h}^*}\right] \leq   \sum_{h\in\mathcal{H}^*}\mathbb{P}\left[\hat{g}_{h}(x)<\frac{g_{\inf}}{2} \right]\leq B_{11}
    n^{-12}
\end{equation}
by Proposition \ref{CT4ProOmega_h^c} and using that $|\mathcal{H}^{*}|<n$.\\

\section{Properties of the estimator $\hat{m}_h$: Bias, Variance, and Oracle Inequality}\label{TC4Control:sesgoYvar}
\vspace{0.5cm}

\subsection{Proof of Proposition \ref{CT4pro:sesgo}: Bias of the estimator $\hat{m}_h$} \label{DemoSesgoHatm}

\begin{eqnarray*}
    \mathbb{E}\{\hat{m}_h (x)\} &=& \mathbb{E}\left\{\frac{1}{n}\sum_{i=1}^n K_h(x-X_i )Y_i  \right\}\\
    &=&
    \mathbb{E}\left\{\frac{1}{n}\sum_{i=1}^n K_h(x-X_i )r(X_i)\right\} + \mathbb{E}\left\{\frac{1}{n}\sum_{i=1}^n K_h(x-X_i )\varepsilon_i  \right\}\\
    &=&
    \mathbb{E}\left\{K_h(x-X)r(X)\right\} + \mathbb{E}\left\{K_h(x-X)\right\} \mathbb{E}\{\varepsilon\}\\
    &=&
    \int_{\mathbb{R}}K_h (x-t)r(t)g(t)dt\\
    &=&
    \int_{\mathbb{R}}K_h (x-t)m(t)dt=
    K_h \ast m(x).
\end{eqnarray*}
Thus, we have
\begin{equation*}
 \mathbb{E}[\hat{m}_h (x)]-m(x)=
   K_h \ast m(x)-m(x)=
   \int_{\mathbb{R}}\left(m(x-hu)-m(x)\right)K\left(u\right)du.
\end{equation*}
From the Taylor-Lagrange formula at order $\mathit{l}_1$ applied to $m(x)$, we have
\begin{equation*}
m(x-hu)=m(x)-hum'(x)+\frac{(-hu)^{2}}{2!}m''(x)-\cdots
+\frac{(-hu)^{\mathit{l}_1}}{\mathit{l}_1 !}m^{(\mathit{l}_1 )}(x-\zeta_u hu)
\end{equation*}
where $\zeta_u\in [0,1]$.

By the above, and since $K$ is a kernel of order $\ord$ with $\ord \geq\mathit{l}_1$, we have
\begin{equation*}
\mathbb{E}[\hat{m}_h (x)]-m(x)=\int_{\mathbb{R}}\frac{(-hu)^{\mathit{l}_1}}{\mathit{l}_1 !}m^{(\mathit{l}_1 )}(x-\zeta_u hu)K(u)du
\end{equation*}

Subtracting the zero term $\int_{\mathbb{R}}\frac{(-uh)^{\mathit{l}_1}}{\mathit{l}_1 !}m^{(\mathit{l}_1 )}(x)K(u)du$ from the previous equation yields

\begin{equation*}
\mathbb{E}[\hat{m}_h (x)]-m(x)=(-1)^{l_1}\frac{h^{\mathit{l}_1}}{\mathit{l}_1 !}\int_{\mathbb{R}}\left(m^{(\mathit{l}_1 )}(x-\zeta_u hu)-m^{(\mathit{l}_1 )}(x)\right)u^{\mathit{l}_1}K(u)du
\end{equation*}

Taking absolute values on both sides and using the hypothesis $m\in\Sigma(\beta_1 ,L_1 )$, we obtain

\begin{align*}
\left|\mathbb{E}[\hat{m}_h (x)]-m(x)\right| &\leq\frac{h^{\beta_1}}{\mathit{l}_1 !}L_1\int_{\mathbb{R}}|u|^{\beta_1}\zeta_u^{\beta_1 -\mathit{l}_1}|K(u)|du\\ &\leq\frac{h^{\beta_1}}{\mathit{l}_1!}L_1\int_{\mathbb{R}}|u|^{\beta_1}|K(u)|du, \quad  \text{since} \quad   \zeta_u^{\beta_1 -\mathit{l}_1}<1\\
&\leq h^{\beta_1}A_{0}.
\end{align*}

\subsection{Proof of Proposition \ref{CT4pro:var}: Variance of the estimator $\hat{m}_h$} \label{DemoVarHatm}

Substituting $Y_i$ by $r(X_i)+\varepsilon_i$, we obtain

\begin{eqnarray} \label{CT4ExpreSigma^2}
\nonumber
     \mathbb{E}\left[\{\hat{m}_h (x)-\mathbb{E}[\hat{m}_h (x)]\}^2\right]
    &=&
    \mathbb{E}\left[\left\{\frac{1}{n}\sum_{i=1}^n K_h(x-X_i )Y_i - \mathbb{E}\left(\frac{1}{n}\sum_{i=1}^n K_h(x-X_i )Y_i \right)\right\}^2\right]\\
    \nonumber
    &=&
    \mathbb{E}\left\{\frac{1}{n}\sum_{i=1}^n \left[ K_h(x-X_i )r(X_i ) - \mathbb{E}(K_h(x-X_i )r(X_i ))\right] \right.\\
    \nonumber
    &+&
    \left. \frac{1}{n}\sum_{i=1}^n K_h(x-X_i )\varepsilon_i \right\}^2\\
    &=&
    I_1 +I_2 +I_3
\end{eqnarray}
where

\begin{itemize}
     \item $A(X_1, \cdots ,X_n )=\sum_{i=1}^n \left[ K_h(x-X_i )r(X_i ) - \mathbb{E}(K_h(x-X_i )r(X_i ))\right]$ 
      \item $I_1 =  \mathbb{E}\left\{\frac{1}{n}\sum_{i=1}^n \left[ K_h(x-X_i )r(X_i ) - \mathbb{E}(K_h(x-X_i )r(X_i ))\right]\right\}^2$
      \item $I_2 = \mathbb{E}\left\{2 \frac{1}{n}A(X_1, \cdots ,X_n ) \frac{1}{n}\sum_{i=1}^n K_h(x-X_i )\varepsilon_i \right\}$
      \item $I_3 =\mathbb{E}\left\{\frac{1}{n}\sum_{i=1}^n K_h(x-X_i )\varepsilon_i \right\}^2$
\end{itemize}

We proceed to bound $I_1$, $I_2$, and $I_3$.

\begin{eqnarray} \label{CT4CotaI3}
    \nonumber
    I_3 &=& \frac{1}{n^2}\mathbb{E}\left\{\sum_{i=1}^n K_h(x-X_i )\varepsilon_i \right\}^2\\
    \nonumber
    &=&
    \frac{1}{n^2}\sum_{i=1}^n \mathbb{E}\left\{K_h(x-X_i )\varepsilon_i \right\}^2 +  \frac{1}{n^2}\sum_{i\neq j}^n \mathbb{E}\left\{K_h(x-X_i )\varepsilon_i K_h(x-X_j )\varepsilon_j \right\}\\
    \nonumber
    &\leq&
    \frac{1}{n^2}\sum_{i=1}^n \frac{B_{1}}{h} +  \frac{1}{n^2}\sum_{i\neq j}^n \mathbb{E}\left\{K_h(x-X_i ) K_h(x-X_j ) \right\}\mathbb{E}( \varepsilon_i\varepsilon_j)\\
    &=&
    \frac{B_{1}}{nh}
\end{eqnarray}
by the first inequality of Lemma \ref{CT4L1}. 

\begin{eqnarray} \label{CT4CotaI2}
   \nonumber
    I_2 &=& \frac{2}{n^2} \mathbb{E}\left\{ A(X_1, \cdots ,X_n )\sum_{i=1}^n K_h(x-X_i )\varepsilon_i \right\}\\
    \nonumber
    &=&
    \frac{2}{n^2}  \sum_{i=1}^n \mathbb{E}\left\{A(X_1, \cdots ,X_n )K_h(x-X_i )\varepsilon_i \right\}\\
    \nonumber
    &=&
    \frac{2}{n^2}  \sum_{i=1}^n \mathbb{E}\left\{A(X_1, \cdots ,X_n )K_h(x-X_i ) \right\}\mathbb{E}\{\varepsilon_i\}=0
\end{eqnarray}

\begin{eqnarray} \label{CT4ExpreI1}
   \nonumber
   I_1 &=&  \frac{1}{n^2}\mathbb{E}\left\{\sum_{i=1}^n \left[ K_h(x-X_i )r(X_i ) - \mathbb{E}(K_h(x-X_i )r(X_i ))\right]\right\}^2\\
   &\leq&
   I_4 + I_5
\end{eqnarray}

where

\begin{itemize}
     \item $I_4 = \left| \frac{1}{n^2}\sum_{i=1}^n Var \left[ K_h(x-X_i )r(X_i )\right] \right|$ 
      \item $I_5 = \left| \frac{1}{n^2}\sum_{i\neq j} Cov \left[K_h(x-X_i )r(X_i ), K_h(x-X_j )r(X_j )\right] \right|$
\end{itemize}

We proceed to bound $I_4$ and $I_5$. By the second inequality of Lemma \ref{CT4L1}, we have
\begin{equation} \label{CT4CotaI4}
   \nonumber
   I_4 \le  \frac{1}{n^2}\sum_{i=1}^n \mathbb{E} \left[ \left( K_h(x-X_i )r(X_i ) \right)^2\right]\le 
   \frac{B_{2}}{nh}. 
\end{equation}

\begin{eqnarray} \label{CT4ExpreI5}
   \nonumber
   I_5 &\leq& \frac{1}{n^2}\sum_{i\neq j}^n \left| Cov \left[K_h(x-X_i )r(X_i ), K_h(x-X_j )r(X_j )\right] \right|\\
   \nonumber
   &=&
   \frac{2}{n^2}\sum_{i=1}^{n-1}\sum_{j=i+1}^n \left| Cov \left[K_h(x-X_i )r(X_i ), K_h(x-X_j )r(X_j )\right] \right|\\
    \nonumber
    &=&
   \frac{2}{n^2}\sum_{i=1}^{n-1}\sum_{k=1}^{n-i} |Cov (\tilde{X}_i, \tilde{X}_{i+k})|\\
    &\le&
   \frac{2}{n^2}\sum_{i=1}^{n-1}\left( \sum_{k=1}^{u_n} |Cov (\tilde{X}_i, \tilde{X}_{i+k})|+\sum_{k=u_n +1}^{{+\infty}} |Cov (\tilde{X}_i, \tilde{X}_{i+k})|\right)
\end{eqnarray}
where $1\leq u_n \leq n-1$ and $\tilde{X}_i =K_h(x-X_i )r(X_i )$, $\forall i$.\\

We analyze the term $\sum_{k=1}^{u_n} |Cov (\tilde{X}_i, \tilde{X}_{i+k})|$ in equation \eqref{CT4ExpreI5}.
\begin{eqnarray} \label{CT4CotaI_i}
   \nonumber
   \sum_{k=1}^{u_n} |Cov (\tilde{X}_i, \tilde{X}_{i+k})| &\leq& \sum_{k=1}^{u_n} |\mathbb{E}[\tilde{X}_i \tilde{X}_{i+k}]|+|\mathbb{E}[\tilde{X}_i]\mathbb{E}[\tilde{X}_{i+k}]|\\
   \nonumber
   &\leq&
   \sum_{k=1}^{u_n} \mathbb{E}[|\tilde{X}_i \tilde{X}_{i+k}|]+\mathbb{E}[|\tilde{X}_i |]\mathbb{E}[|\tilde{X}_{i+k}|]\\
   \nonumber
    &=&
    \sum_{k=1}^{u_n} \left\{ \mathbb{E}[|K_h(x-X_i )r(X_i ) K_h(x-X_{i+k} )r(X_{i+k} )|]+\left( \mathbb{E}[|\tilde{X}_i|] \right)^2 \right\}\\
    \nonumber
     &\leq&
    \sum_{k=1}^{u_n} \left\{ (r_{\sup})^2\mathbb{E}[|K_h(x-X_i ) K_h(x-X_{i+k})|]+\left( r_{\sup}g_{\sup} \|K\|_1 \right)^2 \right\}\\
    \nonumber
     &=&
    \sum_{k=1}^{u_n} \left\{ (r_{\sup})^2Q\|K\|_1^2 + (r_{\sup})^2 (g_{\sup})^2 \|K\|_1^2 \right\}=A_3 u_n
\end{eqnarray}
by the third and fourth inequalities of Lemma \ref{CT4L1}.\\

Using $(H_3)$ and that $ |\tilde{X}_i |\le 
   \frac{r_{\sup}\|K\|_{\infty}}{h}$, we obtain
\begin{eqnarray*}
   |cov(\tilde{X}_i ,\tilde{X}_{i+k})| &\leq& 4 a^k \frac{(r_{\sup})^2\|K\|_{\infty}^2}{h^2}
\end{eqnarray*}

By this last inequality, we have
\begin{eqnarray} \label{CT4CotaJ_i}
   \nonumber
   \sum_{k=u_n +1}^{{+\infty}} |Cov (\tilde{X}_i, \tilde{X}_{i+k})| &\leq& \frac{4(r_{\sup})^2\|K\|_{\infty}^2}{h^2}\sum_{k=u_n +1}^{\infty} a^k\\
   \nonumber
   &=&
   \frac{4(r_{\sup})^2\|K\|_{\infty}^2}{h^2}\frac{a^{u_n +1}}{1-a}\\
   &=&
   A_4\frac{a^{u_n +1}}{h^2}.
\end{eqnarray}

By equation \eqref{CT4ExpreI5} and inequalities \eqref{CT4CotaI_i} and \eqref{CT4CotaJ_i}, we have
\begin{eqnarray*}
   I_5  &\leq&  \frac{2}{n^2}\sum_{i=1}^{n-1}\left( A_3 u_n + A_4\frac{a^{u_n +1}}{h^2}\right)\\
   &\leq&
   \frac{1}{n}\left(2A_3 u_n + 2A_4\frac{a^{u_n +1}}{h^2}\right)
\end{eqnarray*}

The variable $u_n$ in the summations of equation \eqref{CT4ExpreI5} is taken as $u_n = \left[\frac{1}{h|\log a|(\log n)^\frac{1}{2}}\right]$, and using the facts that $[x] \leq x \leq [x]+1$ and $0<a<1$, we continue bounding the previous inequality:

\begin{eqnarray}
   \nonumber
   I_5  &\leq& 
   \frac{1}{n}\left(2A_3 \frac{1}{h|\log a|(\log n)^\frac{1}{2}} + 2A_4\frac{a^{\frac{1}{h|\log a|(\log n)^\frac{1}{2}}}}{h^2}\right)\\
   \nonumber
   &=&
   \frac{1}{nh}\left(2A_3 \frac{1}{|\log a|(\log n)^\frac{1}{2}} + 2A_4\frac{1}{h}a^{\frac{1}{h|\log a|(\log n)^\frac{1}{2}}}\right)
\end{eqnarray}

By Lemma \ref{CT4L2} and the previous inequality, we have that for each $n\geq 4$,

\begin{eqnarray} \label{CT4CotaI5}
    \nonumber
   I_5  &\leq& \frac{1}{nh}\left(\frac{2A_3}{|\log a|} \frac{1}{(\log n)^\frac{1}{2}} + 2A_4 \frac{10}{(\log n)^\frac{5}{2}}\right)\\
    &\leq&
    \left(\frac{2A_3}{|\log a|} + 20A_4\right)\frac{(\log n)^{-\frac{1}{2}}}{nh}.
\end{eqnarray}
By equations \eqref{CT4ExpreI1}, \eqref{CT4CotaI4}, and \eqref{CT4CotaI5}, for all $n\geq 4$,
\begin{equation} \label{CT4CotaI1}
   I_1  \leq
  \frac{B_{2}}{nh}+ \left(\frac{2A_3}{|\log a|} + 20A_4\right)\frac{(\log n)^{-\frac{1}{2}}}{nh}.
\end{equation}
Finally, by inequalities \eqref{CT4CotaI3} and \eqref{CT4CotaI1}, we have
\begin{equation}\label{CT4var:hatm}
   \nonumber
  \mathbb{E}\left[(\hat{m}_h (x)-\mathbb{E}[\hat{m}_h (x))^2\right] = I_1 +I_3\leq 
    \frac{A_{1}}{nh}+A_{2}\frac{(\log n)^{-\frac{1}{2}}}{nh}.
\end{equation}

\subsection{Proof of Theorem \ref{CT4teo:D-Oraculo}: Oracle Inequality for $\hat{m}_{\hat{h}}$} \label{DemoOracuHatm}

{Following the same scheme as the proof of Theorem 2 in Section 4.2 of \citet{mpadrino2024adaptiregre}, the estimator first satisfies}
\begin{eqnarray}\label{CT4CotaErrorCM0}
\nonumber
\lefteqn{\left(\mathbb{E}\left(m(x)-\hat{m}_{\hat{h}}(x)\right)^2\right)^{\frac{1}{2}}}\\
\nonumber
&\leq& \left(\mathbb{E}\left(m(x)-\hat{m}_{h}(x)\right)^2\right)^\frac{1}{2}+2\left(\mathbb{E}(T_{1}^{2})\right)^{\frac{1}{2}}+2\left(\mathbb{E}(T_{2}^{2})\right)^{\frac{1}{2}}+ 2C(h)\|K\|_1 +2V(h)\\
&&
\end{eqnarray}
where
$$T_1=\max_{h'\in\mathcal{H}}\{|\hat{m}_{h'}(x)-\mathbb{E}(\hat{m}_{h'}(x))|-V_1(h')\}_+$$
$$T_2=\max_{h'\in\mathcal{H}}\{|\hat{m}_{h,h'}(x)-\mathbb{E}(\hat{m}_{h,h'}(x))|-V_2(h')\}_+$$
$$C(h)=\max_{u\in B(x)}\left|K_h\ast m(u)-m(u)\right|$$
with 
$$V_{1}(h) = \sqrt{ 2\gamma A_{1}}( 1 + \delta_n ) \frac{(\log n)^\frac{1}{2}}{(nh)^\frac{1}{2}} \quad \text{and} \quad V_2(h) = \sqrt{ 2\gamma A_{1}}\|K\|_{1}( 1 + \delta_n ) \frac{(\log n)^\frac{1}{2}}{(nh)^\frac{1}{2}}.$$

{More especifically Inequality \eqref{CT4CotaErrorCM0} is obtained by following line by line the proof of Proposition 3 in Appendix C.1 of \citet{mpadrino2024adaptiregre}, substituting $\hat{r}_h$ by $\hat{m}_h$ and using that $\mathbb{E}(\hat{m}_h(x))=K_h\ast m(x)$.}\\

Then, $T_1$ and $T_2$ satisfy for $n$ sufficiently large
\begin{equation} \label{CT4CotaET_1^2}
\mathbb{E}(T_1^2) \leq A_5\frac{\log n}{n}
\end{equation}
\begin{equation} \label{CT4CotaET_2^2}
\mathbb{E}(T_2^2) \leq A_5\|K\|_1^2\frac{\log n}{n}
\end{equation}
where 
$$A_5=B_9+2B_{10}.$$
with
{$$B_{9}= \frac{2\|K\|^2_{\infty}}{(2\pi)^{\frac{1}{4}}}\left(2(r_{\sup})^4 +3\sigma^4\right)^{\frac{1}{2}}$$}
$$B_{10} = \frac{2^4 3}{1-e^{-\frac{\gamma-2}{2}}}(A_{1} + B_{8}) + \frac{5!2^{18} 3 B^2}{1-e^{-\frac{\gamma-2}{2}}} (\sigma + r_{\sup})^2$$
{Again, \eqref{CT4CotaET_1^2} and \eqref{CT4CotaET_2^2} are obtained by following line by line the proof of Proposition 4 located in Appendices C.2 and C.3 of \citet{mpadrino2024adaptiregre}} and using Proposition~\ref{CT4TDBernstein}, Corollary~\ref{CT4CDBernstein}, and inequality (ii) of Lemma~\ref{CT4MomenEta1,2y4}.  
\ 
Now, using that for $h\in\mathcal{H}$, $|\mathbb{E}(\hat{m}_h (x)) - m(x)|=|K_h \ast m(x) -m(x)| \leq C(h)$ and Proposition \ref{CT4pro:var}, we have
\begin{eqnarray} \label{CT4ECMTilderh}
    \nonumber
    \left(\mathbb{E}\left[ \left(m(x)-\hat{m}_{h}(x)\right)^2\right]\right)^\frac{1}{2} &\leq&  C(h)+
    \left( \frac{A_{1}}{nh}+A_2\frac{(\log n)^{-\frac{1}{2}}}{nh} \right)^\frac{1}{2} \\
    &\leq&  C(h)+ (A_{1}+A_{2})^{1/2} (nh)^{-1/2}
\end{eqnarray}
\ 
Now, by inequalities \eqref{CT4CotaErrorCM0} and \eqref{CT4ECMTilderh}, we have
\ 
\begin{eqnarray*}
    \nonumber
    \lefteqn{\left(\mathbb{E}\left(m(x)-\hat{m}_{\hat{h}}(x)\right)^2\right)^{\frac{1}{2}}}\\
    &\leq& C(h)+ (A_{1}+A_{2})^{1/2} (nh)^{-1/2} + 2C(h)\|K\|_1 +2 V(h)  +2\left(\mathbb{E}(T_{1}^{2})\right)^{\frac{1}{2}}+2\left(\mathbb{E}(T_{2}^{2})\right)^{\frac{1}{2}}\\
    &\leq&
    A_6 C(h) + A_7 ( 1 + \delta_n ) \frac{(\log n)^{\frac{1}{2}}}{(nh)^{\frac{1}{2}}} +2\left(\mathbb{E}(T_{1}^{2})\right)^{\frac{1}{2}}+2\left(\mathbb{E}(T_{2}^{2})\right)^{\frac{1}{2}}\\
\end{eqnarray*}
Using inequalities \eqref{CT4CotaET_1^2} and \eqref{CT4CotaET_2^2} allows us to obtain 
$$2\left(\mathbb{E}(T_{1}^{2})\right)^{\frac{1}{2}}+2\left(\mathbb{E}(T_{2}^{2})\right)^{\frac{1}{2}}\le  A_8 \frac{(\log n)^\frac{1}{2}}{n^\frac{1}{2}}$$ for a constant $A_{8}=2\sqrt{A_5}(\|K\|_{1}+1)$ depending on the kernel $K$, $r_{sup}$, $Q$, $g_{sup}$, $\sigma$, $a$, and $\gamma$, which completes the proof.

\section{Inequalities for the regression estimator and adaptivity of $\hat{r}_{\hat{h}/\hat{h}^*}$}\label{TC4UltimaParte}
\vspace{0.5cm}

\subsection{Proof of Proposition \ref{CT4DesConsVel}} \label{DemoPropo5}

Let $n$ be such that $\frac{1}{\log n} < \frac{g_{\inf}}{2}$.  {For $h\in[h_{min},h_{max}]$ and $h^*\in[h_{min^*},h_{max}^*]$,} we have

$$\left(r(x)-\hat{r}_{h/h^*}(x)\right)^2 \leq \frac{2}{(\bar{g}_{h^*}(x))^2}(m(x)-\hat{m}_{h}(x))^2+2\left(\frac{r(x)}{\bar{g}_{h^*}(x)}\right)^2(\bar{g}_{h^*}(x)-g(x))^2$$
and therefore,

\begin{eqnarray} \label{CT4DesConsVel1}
    \nonumber
    \E\left[\left(r(x)-\hat{r}_{h/h^*}(x)\right)^2\right] &\leq& 2\E\left[\frac{1}{(\bar{g}_{h^*}(x))^2}(m(x)-\hat{m}_{h}(x))^2\right]\\
    &+&
    2(r_{\sup})^2\E\left[\frac{1}{(\bar{g}_{h^*}(x))^2}(\bar{g}_{h^*}(x)-g(x))^2\right]
\end{eqnarray}

Recall that $\Bar{g}_{h^*}(x)=\max\left\{\hat{g}_{h^*}(x),\frac{1}{\log n}\right\}$. Consider the event $\Omega_{h^*}=\left\{\hat{g}_{h^*}(x) \geq \frac{g_{\inf}}{2}\right\}$. 

On $\Omega_{h^*}$, we have $\bar{g}_{h^*}(x)=\hat{g}_{h^*}(x)$ and
\begin{equation} \label{CT4Cota2(1/Barg)}
    \frac{1}{\bar{g}_{h^*}(x)} \leq \frac{2}{g_{\inf}}.
\end{equation}

From inequality \eqref{CT4DesConsVel1}, we analyze the first term:

\begin{eqnarray} \label{CT4DesConsVel2}
    \nonumber
    \E\left[\frac{1}{(\bar{g}_{h^*}(x))^2}(m(x)-\hat{m}_{h}(x))^2\right] &=& \E\left[\frac{1}{(\bar{g}_{h^*}(x))^2}(m(x)-\hat{m}_{h}(x))^2\mathbf{1}_{\Omega_{h^*}}\right]\\
    &+&
    \E\left[\frac{1}{(\bar{g}_{h^*}(x))^2}(m(x)-\hat{m}_{h}(x))^2\mathbf{1}_{\Omega_{h^*}^c}\right]
\end{eqnarray}

By inequality \eqref{CT4Cota2(1/Barg)} and the first term of equation \eqref{CT4DesConsVel2}, we have

\begin{eqnarray} \label{CT4DesConsVel3}
    \nonumber
    \E\left[\frac{1}{(\bar{g}_{h^*}(x))^2}(m(x)-\hat{m}_{h}(x))^2\mathbf{1}_{\Omega_{h^*}}\right] &\leq& \frac{4}{(g_{\inf})^2}\E\left[(m(x)-\hat{m}_{h}(x))^2\mathbf{1}_{\Omega_{h^*}}\right]\\
    &\leq&
    \frac{4}{(g_{\inf})^2}\E\left[(m(x)-\hat{m}_{h}(x))^2\right]
\end{eqnarray}

By definition of $\bar{g}_{h^*}$ and the second term of equation \eqref{CT4DesConsVel2}, we have

\begin{eqnarray} \label{CT4DesConsVel4}
    \nonumber
    \lefteqn{
    \E\left[\frac{1}{(\bar{g}_{h^*}(x))^2}(m(x)-\hat{m}_{h}(x))^2\mathbf{1}_{\Omega_{h^*}^c}\right]}\\
    \nonumber
    &\leq& (\log n)^2\E\left[(m(x)-\hat{m}_{h}(x))^2\mathbf{1}_{\Omega_{h^*}^c}\right]\\
    \nonumber
    &\leq&
    (\log n)^2\E\left[\left(2(m(x))^2 + 2(\hat{m}_{h}(x))^2\right)\mathbf{1}_{\Omega_{h^*}^c}\right]\\
    \nonumber
    &=&
    2(\log n)^2\left\{ \E\left[(m(x))^2\mathbf{1}_{\Omega_{h^*}^c}\right] + \E\left[(\hat{m}_{h}(x))^2\mathbf{1}_{\Omega_{h^*}^c}\right] \right\}\\
    \nonumber
    &\leq&
    2(\log n)^2\left\{ (r_{\sup}g_{\sup})^2\mathbb{P}\left[\Omega_{h^*}^c \right] +
   I_1\right\}\\
   &\le & 2(r_{\sup}g_{\sup})^{2}B_{11} n^{-4}+2(\log n)^2I_1
\end{eqnarray}
where the last line is a consequence of Proposition \ref{CT4ProOmega_h^c}, and $I_{1}=\E\left[\left( \hat{m}_{h}(x) \right)^2\mathbf{1}_{\Omega_{h^*}^c}\right]$.\\

Now, using Lemma~\ref{CT4MomenY^2,4,8} we deduce that
{
\begin{equation*}
\begin{aligned}
\mathbb E[\widehat m_h(x)^4]
&\leq \frac{\|K\|_\infty^4}{h_{\min}^4}
 \mathbb E\!\left[\left(\frac1n\sum_{i=1}^n
  |Y_i|\mathbf1_{\{|X_i-x|\leq h_{\max}\}}\right)^4\right]\\
&\leq \frac{\|K\|_\infty^4}{h_{\min}^4}
 2(2r_{\sup}^4+3\sigma^4),
\end{aligned}
\end{equation*}
which implies that:
\begin{equation}\label{CT4DesConsVel7}
\begin{aligned}
I_1&\leq \bigl(\mathbb E[\widehat m_h(x)^4]\bigr)^{1/2}
          \mathbb P(\Omega_{h^*}^c)^{1/2}\\
&\leq
\|K\|_\infty^2\sqrt{2(2r_{\sup}^4+3\sigma^4)B_{11}}\,
 h_{\min}^{-2}n^{-13/2},
\end{aligned}
\end{equation}
}
where the last line is a consequence of Proposition \ref{CT4ProOmega_h^c}, and $h_{min}=\frac{(\log n)^8}{n}\geq \frac{\log n}{n}$. Thus,
\begin{equation} \label{CT4DesConsVel8}
    \E\left[\frac{1}{(\bar{g}_{h^*}(x))^2}(m(x)-\hat{m}_{h}(x))^2\mathbf{1}_{\Omega_{h^*}^c}\right]\le \Bar{F}_{1} n^{-4},
\end{equation}
 where $\Bar{F}_{1}=  2B_{11}(r_{sup}g_{sup})^2 + 2\|K\|^2_{\infty}B_{11}^{1/2}(2(2r_{sup}^4+3\sigma^4))^{1/2}$.

A bound for the first term on the right-hand side of inequality \eqref{CT4DesConsVel1} is given using inequalities \eqref{CT4DesConsVel2}, \eqref{CT4DesConsVel3}, and \eqref{CT4DesConsVel8}:
\begin{equation} \label{CT4DesConsVel9}
    \E\left[\frac{1}{(\bar{g}_{h^*}(x))^2}(m(x)-\hat{m}_{h}(x))^2\right] \leq \frac{4}{(g_{inf})^2}\E\left[(m(x)-\hat{m}_{h}(x))^2\right] + \frac{\Bar{F}_1}{n^4}.
\end{equation}

Now, a bound for the second term on the right-hand side of inequality \eqref{CT4DesConsVel1} is given:

\begin{eqnarray} \label{CT4DesConsVel10}
    \nonumber
    \E\left[\frac{1}{(\bar{g}_{h^*}(x))^2}(\Bar{g}_{h^*}(x)-g(x))^2\right] &=& \E\left[\frac{1}{(\bar{g}_{h^*}(x))^2}(\Bar{g}_{h^*}(x)-g(x))^2\mathbf{1}_{\Omega_{h^*}}\right]\\
    &+&
    \E\left[\frac{1}{(\bar{g}_{h^*}(x))^2}(\Bar{g}_{h^*}(x)-g(x))^2\mathbf{1}_{\Omega_{h^*}^c}\right]
\end{eqnarray}
For the first term on the right-hand side of equation \eqref{CT4DesConsVel10}, using equation \eqref{CT4Cota2(1/Barg)} and the fact that on $\Omega_{h^*}$, $\bar{g}_{h^*}(x)=\hat{g}_{h^*}(x)$, we have
\begin{eqnarray} \label{CT4DesConsVel11}
    \nonumber
    \E\left[\frac{1}{(\bar{g}_{h^*}(x))^2}(\Bar{g}_{h^*}(x)-g(x))^2\mathbf{1}_{\Omega_{h^*}}\right] &\leq& \E\left[\left( \frac{2}{g_{inf}} \right)^{2}(\hat{g}_{h^*}(x)-g(x))^2\mathbf{1}_{\Omega_{h^*}}\right]\\
    &\leq&
    \frac{4}{(g_{inf})^2}\E\left[(\hat{g}_{h^*}(x)-g(x))^2\right]
\end{eqnarray}
For the second term on the right-hand side of equation \eqref{CT4DesConsVel10}, we have
\begin{eqnarray} \label{CT4DesConsVel12}
    \nonumber
    \E\left[\frac{1}{(\bar{g}_{h^*}(x))^2}(\Bar{g}_{h^*}(x)-g(x))^2\mathbf{1}_{\Omega^c_{h^*}}\right] &\leq& (\log n)^2\E\left[(\bar{g}_{h^*}(x)-g(x))^2\mathbf{1}_{\Omega_{h^*}^c}\right]\\
    \nonumber
    &\leq&
    (\log n)^2\E\left[\left(2(\bar{g}_{h^*}(x))^{2}+2(g(x))^2\right)\mathbf{1}_{\Omega_{h^*}^c}\right]\\
    &\leq&
    2(\log n)^2\{J_{1}+J_{2}\}
\end{eqnarray}
where $J_{1}=\E\left[(\bar{g}_{h^*}(x))^{2}\mathbf{1}_{\Omega_{h^*}^c}\right]$ and $J_{2}=\E\left[(g(x))^2\mathbf{1}_{\Omega_{h^*}^c}\right]$.\\

Before bounding the term $J_1$ of inequality \eqref{CT4DesConsVel12}, recall that $\Bar{g}_{h^*}(x)=\max\left\{\hat{g}_{h^*}(x),\frac{1}{\log n}\right\}$.
Using $$|\hat{g}_{h^*}(x)| = \left| \frac{1}{nh^*}\sum_{i=1}^{n} K\left(\frac{x-X_i}{h^*}\right) \right| \leq \frac{\|K\|_{\infty}}{h_{\min}^*}$$
we obtain $\Bar{g}_{h^*}(x) \leq \left( \frac{\|K\|_{\infty}}{h_{\min}^*} \vee \frac{1}{\log n} \right)$. Substituting this into term $J_1$ of inequality \eqref{CT4DesConsVel12}, we have
\begin{equation}  \label{CT4DesConsVel13}
    J_{1} \leq
    \left( \frac{\|K\|^2_{\infty}}{(h_{\min}^*)^2} \vee \frac{1}{(\log n)^2} \right)\mathbb{P}\left[\Omega_{h^*}^c\right]
\end{equation}

Now we proceed to bound the term $J_2$ of inequality \eqref{CT4DesConsVel12}:

\begin{equation} \label{CT4DesConsVel14}
    J_{2} \leq \E\left[(g_{\sup})^2\mathbf{1}_{\Omega_{h^*}^c}\right] = (g_{\sup})^2\mathbb{P}\left[ \Omega_{h^*}^c \right]
\end{equation}

By inequalities \eqref{CT4DesConsVel12}, \eqref{CT4DesConsVel13}, and \eqref{CT4DesConsVel14}, in addition to Proposition \ref{CT4ProOmega_h^c} and since $h^*_{\min}=n^{-1}e^{(\sqrt{\log n}-1)}$, we have

\begin{eqnarray} \label{CT4DesConsVel15}
    \nonumber
    \lefteqn{\E\left[\frac{1}{(\bar{g}_{h^*}(x))^2}(\Bar{g}_{h^*}(x)-g(x))^2\mathbf{1}_{\Omega^c_{h^*}}\right]}\\
    \nonumber
    &\leq& 2(\log n)^2\left\{\left( \frac{\|K\|^2_{\infty}}{(h_{\min}^*)^2} \vee \frac{1}{(\log n)^2} \right) + (g_{\sup})^2\right\}\mathbb{P}\left[ \Omega_{h^*}^c \right]\\
    \nonumber
    &\leq&
    \frac{2B_{11}(\log n)^2}{n^{13}}\left\{\left( \frac{\|K\|^2_{\infty}}{\left( n^{-1}e^{(\sqrt{\log n}-1)} \right)^2} \vee \frac{1}{(\log n)^2} \right) + (g_{\sup})^2\right\}\nonumber\\
    &\leq&
    \frac{\Bar{F}_{2}}{n^{4}}
\end{eqnarray}
where $\Bar{F}_{2}=2B_{11}\left\{\left(\|K\|^2_{\infty} \vee 1\right) + (g_{\sup})^{2}\right\}$.\\

By equations \eqref{CT4DesConsVel10}, \eqref{CT4DesConsVel11}, and \eqref{CT4DesConsVel15}, we have

\begin{equation} \label{CT4DesConsVel16}
    \E\left[\frac{1}{(\bar{g}_{h^*}(x))^2}(\Bar{g}_{h^*}(x)-g(x))^2\right] \leq \frac{4}{(g_{inf})^2}\E\left[(\hat{g}_{h^*}(x)-g(x))^2\right] + \frac{\Bar{F}_{2}}{n^{4}}.
\end{equation}
Finally, by inequalities \eqref{CT4DesConsVel1}, \eqref{CT4DesConsVel9}, and \eqref{CT4DesConsVel16}:

\begin{equation} \label{CT4DesConsVel17}
    \E\left[\left(r(x)-\hat{r}_{h/h^*}(x)\right)^2\right] \leq
    F_{1}\E\left[(m(x)-\hat{m}_{h}(x))^2\right] + F_{2}\E\left[(\hat{g}_{h^*}(x)-g(x))^2\right]+
    \frac{F_{3}}{n^{4}}
\end{equation}
where $F_{1}=\frac{8}{(g_{inf})^2}$, $F_{2}=\frac{8(r_{sup})^2}{(g_{inf})^2}$, and $F_{3}=2\Bar{F}_{1} + 2(r_{sup})^2\Bar{F}_{2}$.

\subsection{Proof of Theorem \ref{CT4adapHatr}} \label{DemoCT4adapHatr}

Similarly to Proposition \ref{CT4DesConsVel}, we have
$$\E\left[\left(r(x)-\hat{r}_{\hat{h}/\hat{h}^*}(x)\right)^2\right]\leq F_1\E\left[\left(m(x)-\hat{m}_{\hat{h}}(x)\right)^2\right]+F_2\E\left[\left(g(x)-\hat{g}_{\hat{h}^*}(x)\right)^2\right]+\frac{F_3}{n^4}$$
using $\mathbb{P}(\Omega^{c}_{\hat{h}^*}) \leq \frac{1}{n^{12}}$ (see Proposition \ref{CT4ProOmega_Hath^c}). 

By Theorem \ref{CT4AdapHatm} and Theorem 2 in \cite{bertin2017pointwise}, we have
\begin{eqnarray*}
\lefteqn{\E\left[\left(\hat{r}_{\hat{h}/\hat{h}^*}(x)-r(x)\right)^2\right]}\\
&\leq& F_{1}\left\{ C^* \left( \frac{n}{\log n} \right)^{-\frac{\beta_1}{2\beta_{1}+1}} \right\}^{2}+F_{2} \left\{ C^{**}\left( \frac{\log n}{n} \right)^{\frac{\beta_2}{2\beta_{2}+1}} \right\}^{2} + \frac{F_3}{n^4}\\
&\leq&{\mathfrak{B}^2} \left( \frac{n}{\log n}\right)^{-\frac{2(\beta_{1}{\wedge}\beta_{2})}{2(\beta_{1}{\wedge}\beta_{2})+1}}
\end{eqnarray*}
where {$\mathfrak{B}=(F_{1}\left({C^*})^2 +F_{2}(C^{**})^{2} +F_{3}\right)^{1/2}$ with $C^{**}$ }a positive constant depending on the kernel $K$, $g_{\sup}$, $Q$, $\beta_2$, $L_2$, and $|\log(a)|$.\\

\section{Calibration}\label{Calibra}
\vspace{0.5cm}

In this appendix, we present complementary numerical results related to the calibration of the parameters $\gamma$ and $\gamma^{*}$ involved in the Goldenshluger - Lepski selection procedures used to estimate $m$ and $g$, respectively.

The choice of the calibration intervals considered in the numerical study was motivated by the empirical observation of a transition between overfitting and underfitting regimes as these parameters vary. This behavior is reminiscent of the phase-transition phenomena associated with the concept of minimal penalty established for the Goldenshluger--Lepski procedure by Lacour and Massart \cite{lacour2016minimalpen}.

Although a theoretical analysis of these phenomena in our setting is beyond the scope of the present work, the results reported below provide empirical evidence of analogous behavior. In particular, the figures in this appendix illustrate how variations in $\gamma$ and $\gamma^{*}$ affect the selected bandwidths and, consequently, the final estimates produced by the adaptive procedure.

It is worth emphasizing that, unlike the classical setting, our framework involves two Goldenshluger--Lepski procedures operating simultaneously. Determining whether critical constants associated with $\gamma$ and $\gamma^{*}$ exist, as well as understanding the interaction between both selection mechanisms in the construction of the adaptive estimator
$\widehat r=\frac{\widehat m}{\widehat g}$,
constitutes an open problem of theoretical interest.

\subsection{Calibration of the Method.} \label{CCali} 

When calibrating the method for a sample $\{(X_{i},Y_{i})\}_{i=1}^{n+2q}$, the first $n$ observations are used to construct the estimators of $g$ and $m$, namely $\hat{g}_{h^*}$ and $\hat{m}_h$, with $h^{*},h\in\mathcal{H}$. The last $q$ observations of the sample, $\{(X_{i},Y_{i})\}_{i=n+q+1}^{n+2q}$, are used to construct a random grid on the interval $[-2,2]$ by rearranging the observations in increasing order with respect to their first component. That is,
$\{(X_{i},Y_{i})\}_{i=n+q+1}^{n+2q}=\{(\Tilde{x}_{i},\Tilde{y}_{i})\}_{i=1}^q$,
where $\Tilde{x}_{i}<\Tilde{x}_{j}$ whenever $i<j$. The calibration of the method is carried out using the grid determined by $\{\Tilde{x}_{i}\}_{i=1}^q$ on the interval $[-2,2]$.\\

The estimator of the regression function $r$ is computed at each point of the grid $\{\Tilde{x}_i\}_{i=1}^{q}$ on the interval $[-2,2]$. For each $\Tilde{x}_i$, two optimal bandwidths $\hat{h}_i$ and $\hat{h}^*_i$ are selected from the collection of bandwidths $\mathcal{H}$. The bandwidth $\hat{h}_{i}$ is obtained according to equation \eqref{Chxi}, where $V(h,x_{i})$ depends on the parameter $\gamma$, while the bandwidth $\hat{h}^*_{i}$ is obtained according to equation \eqref{Ch*xi}, where $V^{*}(h,x_{i})$ depends on the parameter $\gamma^*$. Consequently, the optimal bandwidths $\hat{h}_i$ and $\hat{h}^*_i$ depend on the parameters $\gamma$ and $\gamma^*$, respectively, and this dependence is denoted by $\hat{h}^{\gamma}_i$ and $\hat{h}^{*\gamma^*}_i$.

In practice, to calibrate the method, $\gamma$ and $\gamma^*$ are chosen from equally spaced grids over a subinterval $I$ of $]0,1/10[$. The interval $I$ is selected so that, as the values of $\gamma$ and $\gamma^*$ vary within $I$, the curves described by the estimators $\left\{\hat{m}_{\hat{h}_{i}^{\gamma}}(x_{i})\right\}_{i=1}^q$ and $\left\{\hat{g}_{\hat{h}_{i}^{*\gamma^*}}(x_{i})\right\}_{i=1}^q$ transition from highly irregular curves to smooth curves.\\

The calibration procedure is carried out in three stages. First, an appropriate value of $\gamma^*$ is selected within $I$ so that $\left\{\hat{g}_{\hat{h}_i^{*\gamma^*}}(x_i)\right\}_{i=1}^q$ provides a good estimate of $\{g(x_i)\}_{i=1}^q$. Second, an appropriate value of $\gamma$ is selected within $I$ so that $\left\{\hat{m}_{\hat{h}_i^{\gamma}}(x_i)\right\}_{i=1}^q$ provides a good estimate of $\{m(x_i)\}_{i=1}^q$. Finally, using the results obtained in Stages 1 and 2, the estimator of the regression function $r$ is constructed.\\

\textbf{Stage 1 (Selection of $\gamma^*$):} Using the \texttt{hist} function in R, a histogram is constructed from the data $\{X_i\}_{i=1}^n$, and the corresponding histogram function $H_{g}(x)$ is defined for each $x\in [-2,2]$. Next, $Error_{g}(\gamma^*) = \sum_{i=1}^q d_{i}\left(\hat{g}_{\hat{h}_i^{*\gamma^*}}(\tilde{x}_{i})-H_{g}(\tilde{x}_{i})\right)^2$ is computed for each $\gamma^*\in I$, where
$d_{1}=\frac{\tilde{x}_{1}+\tilde{x}_{2}}{2}+2$,
$d_{q}=2-\frac{\tilde{x}_{q-1}+\tilde{x}_{q}}{2}$,
and
$d_{i}=\frac{\tilde{x}_{i+1}-\tilde{x}_{i-1}}{2}$,
for $i\in\{2,\dotsc,q-1\}$.
Finally, the first calibration stage is completed by taking $\gamma_{g}=\arg\min_{\gamma^*\in I}\left\{Error_{g}(\gamma^{*})\right\}$.\\

\textbf{Stage 2 (Selection of $\gamma$):} Since the model under consideration is
$Y_i=\frac{m(X_{i})}{g(X_{i})}+\varepsilon_i,\quad i=1,\ldots,n$,
it follows that $m(X_{i})\approx Y_{i}g(X_{i})$. Therefore, this calibration stage uses the last $q$ reordered observations of the sample,
$\{(\Tilde{x}_{i},\Tilde{y}_{i})\}_{i=1}^q$,
together with the estimator
$\left\{\hat{g}_{\hat{h}_i^{*\gamma_g}}(\Tilde{x}_i)\right\}_{i=1}^q$
obtained in the previous stage. The quantity $Error_{m}(\gamma) =
\sum_{i=1}^q d_{i}
\left(
\hat{m}_{\hat{h}_i^{\gamma}}(\tilde{x}_{i})
-\Tilde{y}_{i}\hat{g}_{\hat{h}_i^{*\gamma_g}}(\Tilde{x}_i)
\right)^2$
is computed for each $\gamma\in I$, where $\{d_i\}_{i=1}^q$ is defined as in Stage 1. Finally, this calibration stage is completed by taking $\gamma_{m}=\arg\min_{\gamma\in I}\left\{Error_{m}(\gamma)\right\}.$\\

\textbf{Stage 3 (Estimator of $r$):} For the estimation of the regression function $r=\frac{m}{g}$ on the equally spaced grid $\{x_i\}_{i=1}^{s}$ over the interval $[-1,1]$, where
$x_{i}=-1+2\left(\frac{i-1}{s-1}\right)$
for $i=1,\cdots,s$ with $s=21$, the estimator
$\left\{\hat{r}_{\hat{h}_{i}^{\gamma_m}/\hat{h}^{*\gamma_g}_i}(x_{i})\right\}_{i=1}^s$
is considered. To simplify the notation, whenever the GL estimator is said to be calibrated, the dependence of the optimal bandwidths on $\gamma_m$ and $\gamma_g$ is omitted, and the calibrated GL estimator is simply written as
$\left\{\hat{r}_{\hat{h}_{i}/\hat{h}^{*}_i}(x_{i})\right\}_{i=1}^s$.\\

In practice, the interval $I=]0.00000005,0.08[$ was used. Moreover, the equally spaced grid $\gamma_i =0.00000005+(i-1)\frac{0.08-0.00000005}{20}, \quad i=1,\cdots,21$, was considered on this interval. 

\begin{figure}[htp!]
    \begin{center}
        \includegraphics[width=7.66cm]{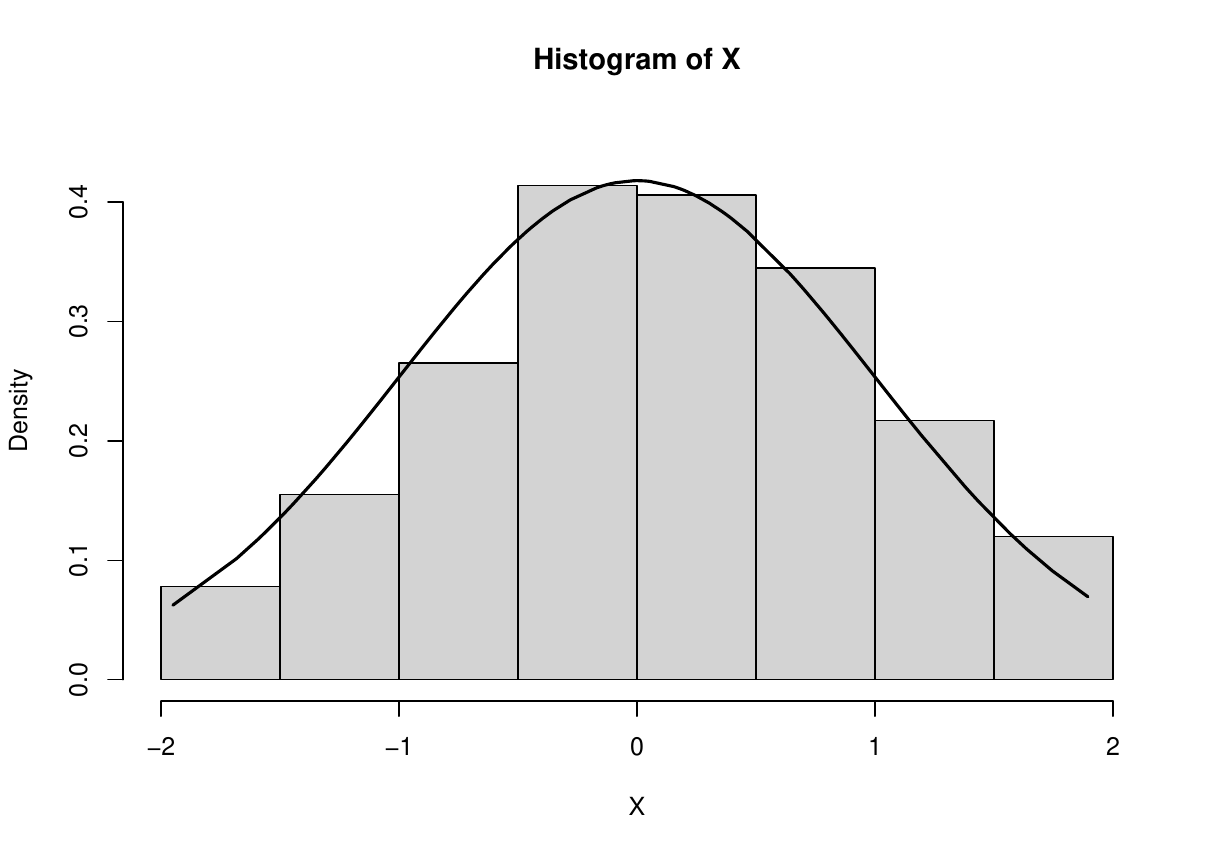}
    \end{center}
   \caption{Histogram constructed from the sample $\{X_i\}_{i=1}^n$ with $n=2000$, together with the density function $g$ evaluated on the random grid $\{\Tilde{x}_i\}_{i=1}^q$, with $q=100$.}
    \label{CHistyg}
\end{figure}

\begin{figure}[htp!] 
	\centering 
	\begin{subfigure}
	\centering
	\includegraphics[width=7.66cm]{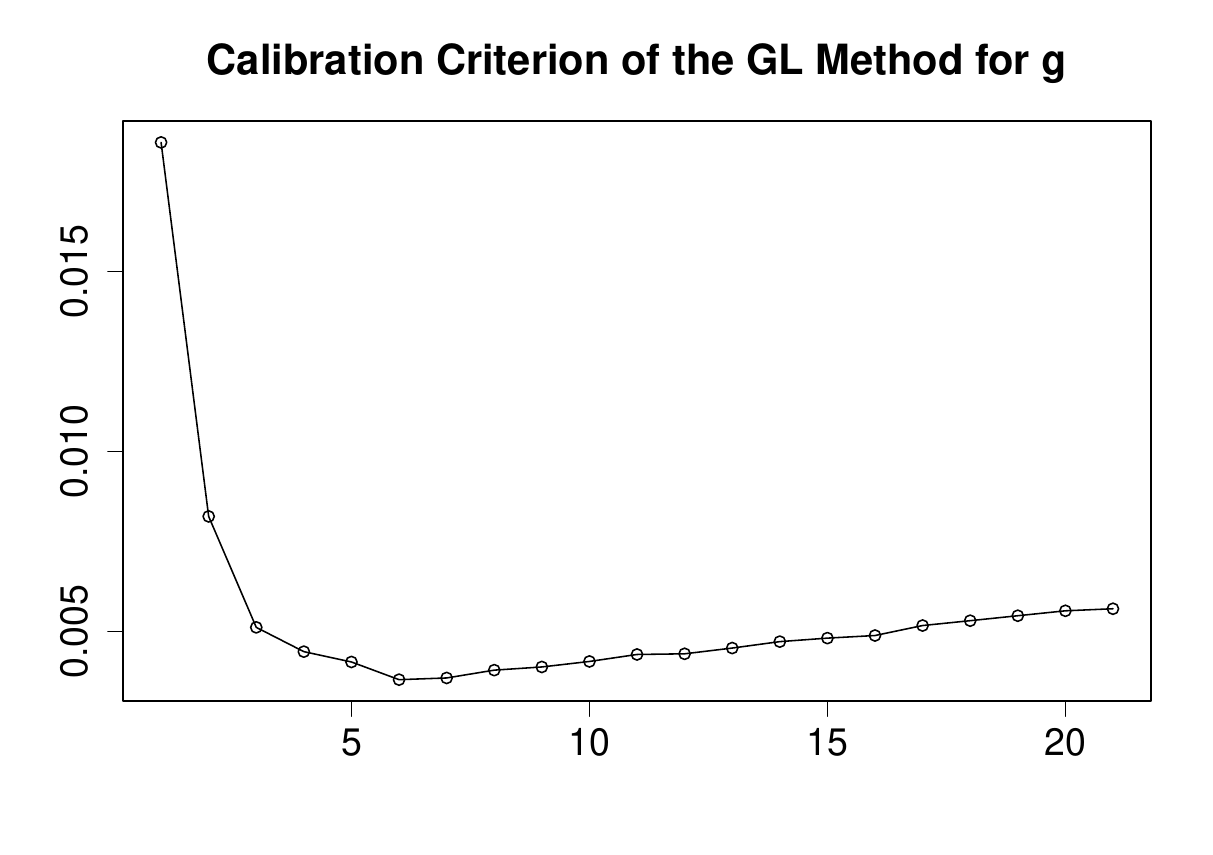}  
	\end{subfigure}
	\begin{subfigure}
	\centering
	\includegraphics[width=7.66cm]{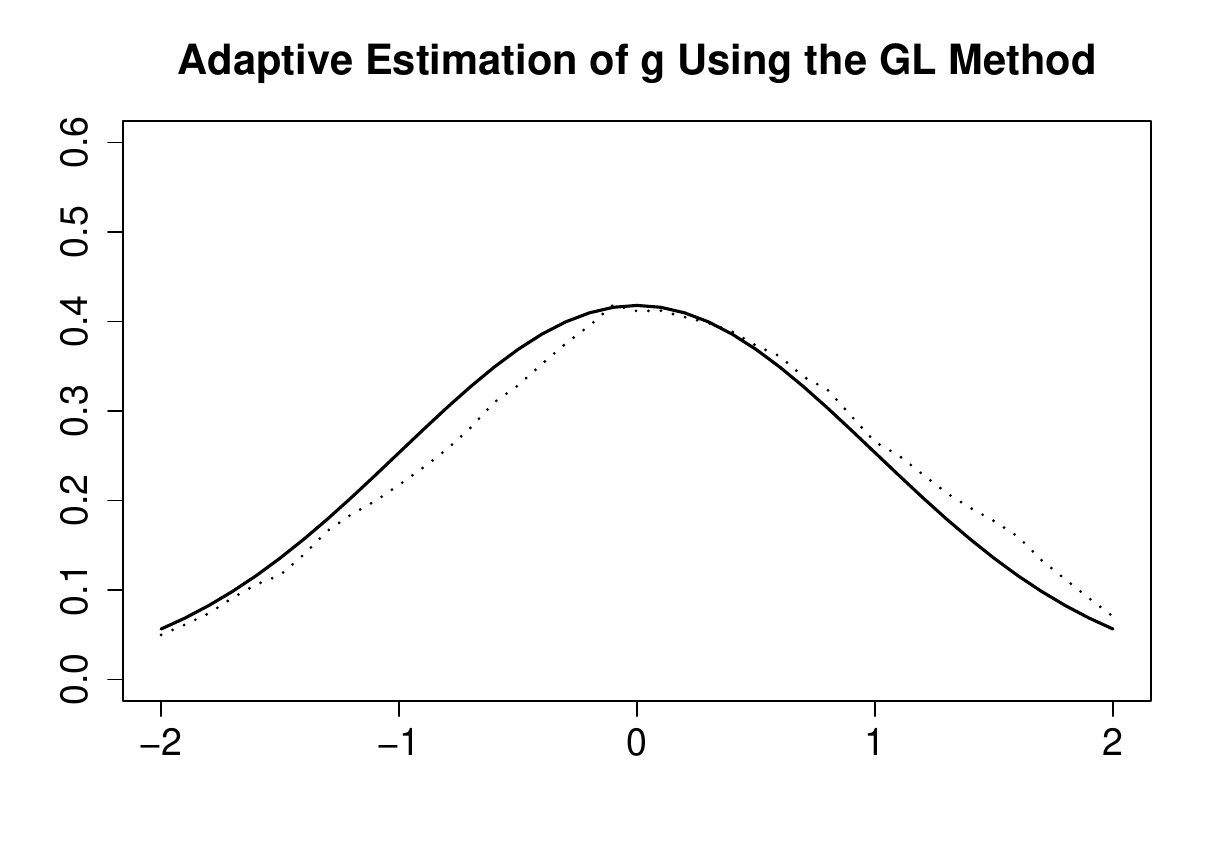}  
	\end{subfigure}
	\caption{Left panel: plot of $\{Error_{g}(\gamma_i^{*})\}_{i=1}^{21}$. Right panel: the density function $g$ is represented by a solid line, while the estimator $\{\hat{g}_{\hat{h}_i^*}(x_i)\}_{i=1}^{41}$ on the interval $[-2,2]$ is represented by dots. The estimator is calibrated at $\gamma_6^*$, corresponding to the value at which the curve in the left panel attains its minimum.}
    \label{CCalibraciong}
\end{figure}

\begin{figure}[htp!] 
	\centering 
	\begin{subfigure}
	\centering
	\includegraphics[width=7.66cm]{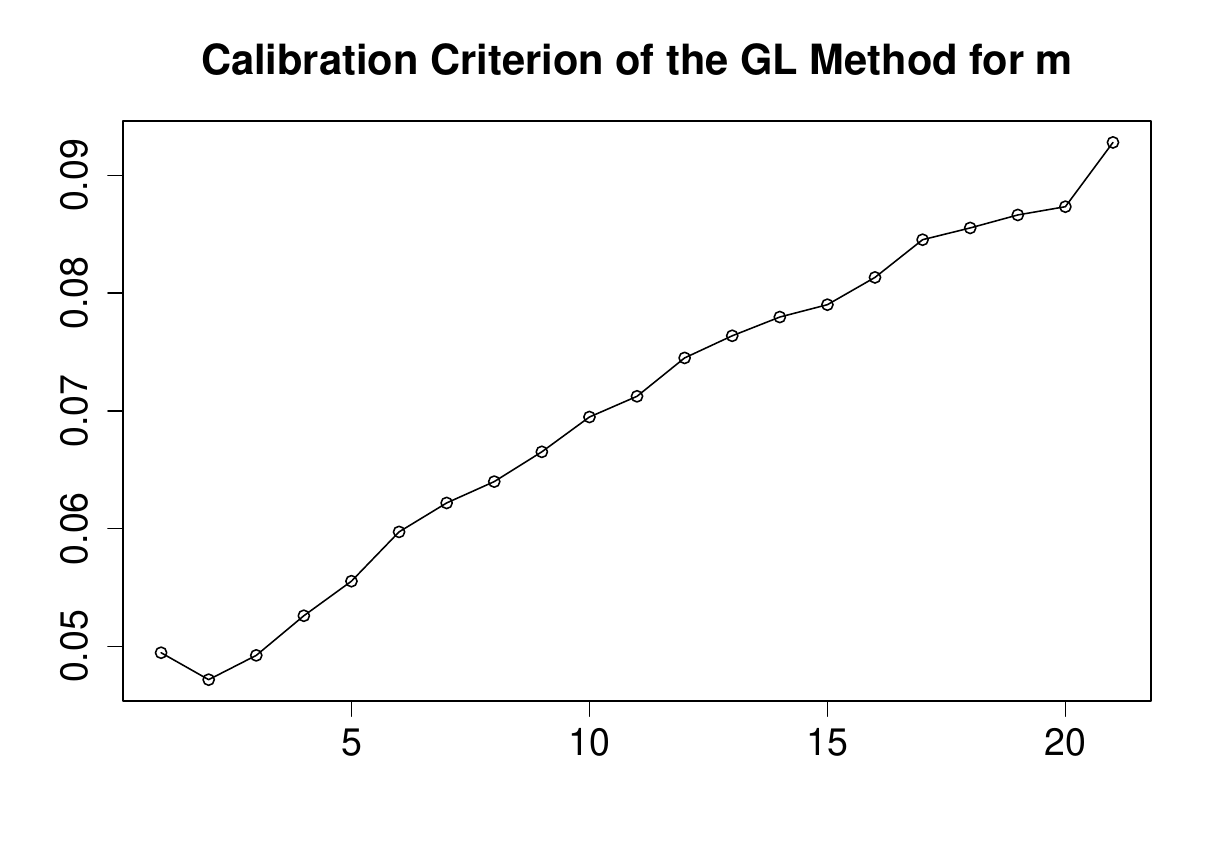}  
	\end{subfigure}
	\begin{subfigure}
	\centering
	\includegraphics[width=7.66cm]{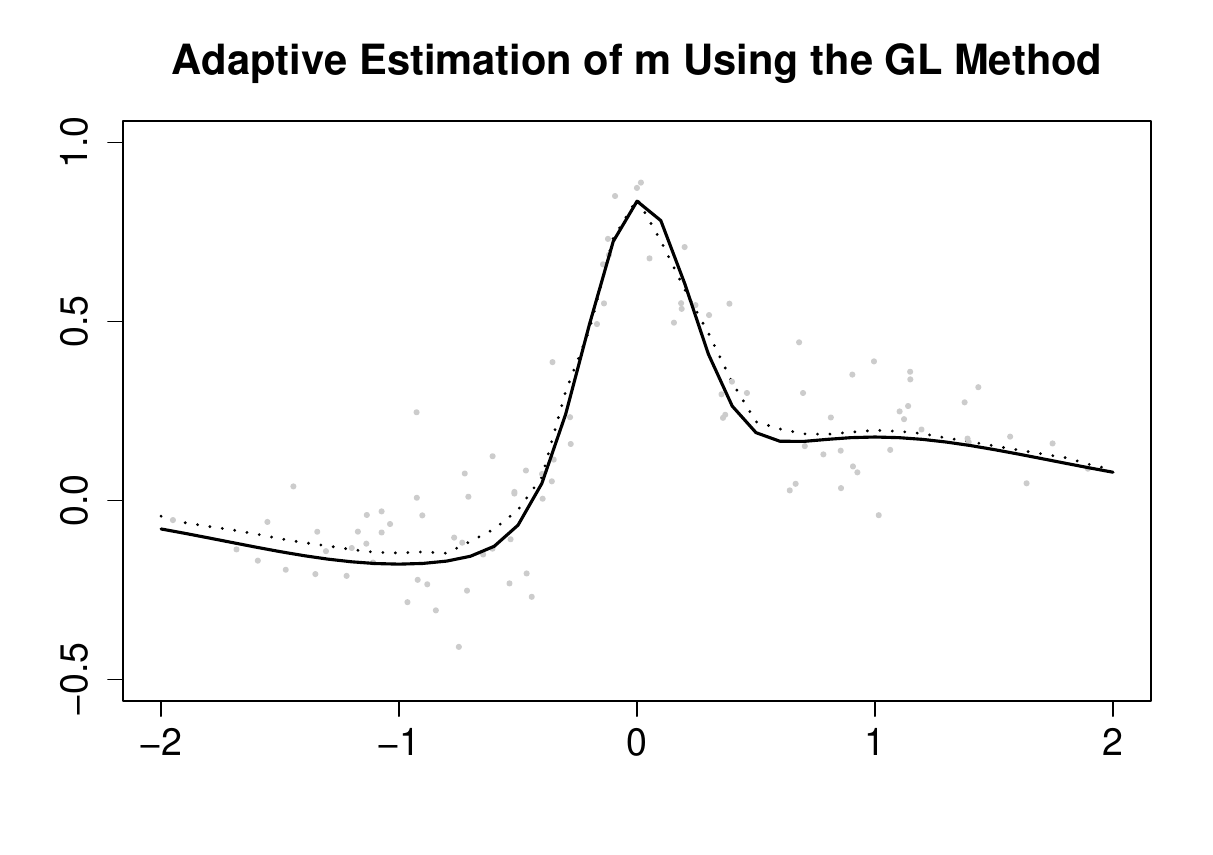}  
	\end{subfigure}
	\caption{Left panel: plot of $\{Error_{m}(\gamma_i)\}_{i=1}^{21}$. Right panel: the point cloud
    $\left\{ \left(\Tilde{x}_i ,\Tilde{y}_{i}\hat{g}_{\hat{h}_i^{*\gamma^*_6}}(\Tilde{x}_i)\right) \right\}_{i=1}^q$
    is displayed. The function $m$ over the interval $[-2,2]$ is represented by a solid line, while the estimator $\{\hat{m}_{\hat{h}_i}(x_i)\}_{i=1}^{41}$ is represented by dots. The estimator is calibrated at $\gamma_2$, corresponding to the value at which the curve in the left panel attains its minimum.}
    \label{CCalibracionm}
\end{figure}

\begin{figure}[htp!] 
	\centering 
	\begin{subfigure}
	\centering
	\includegraphics[width=7.66cm]{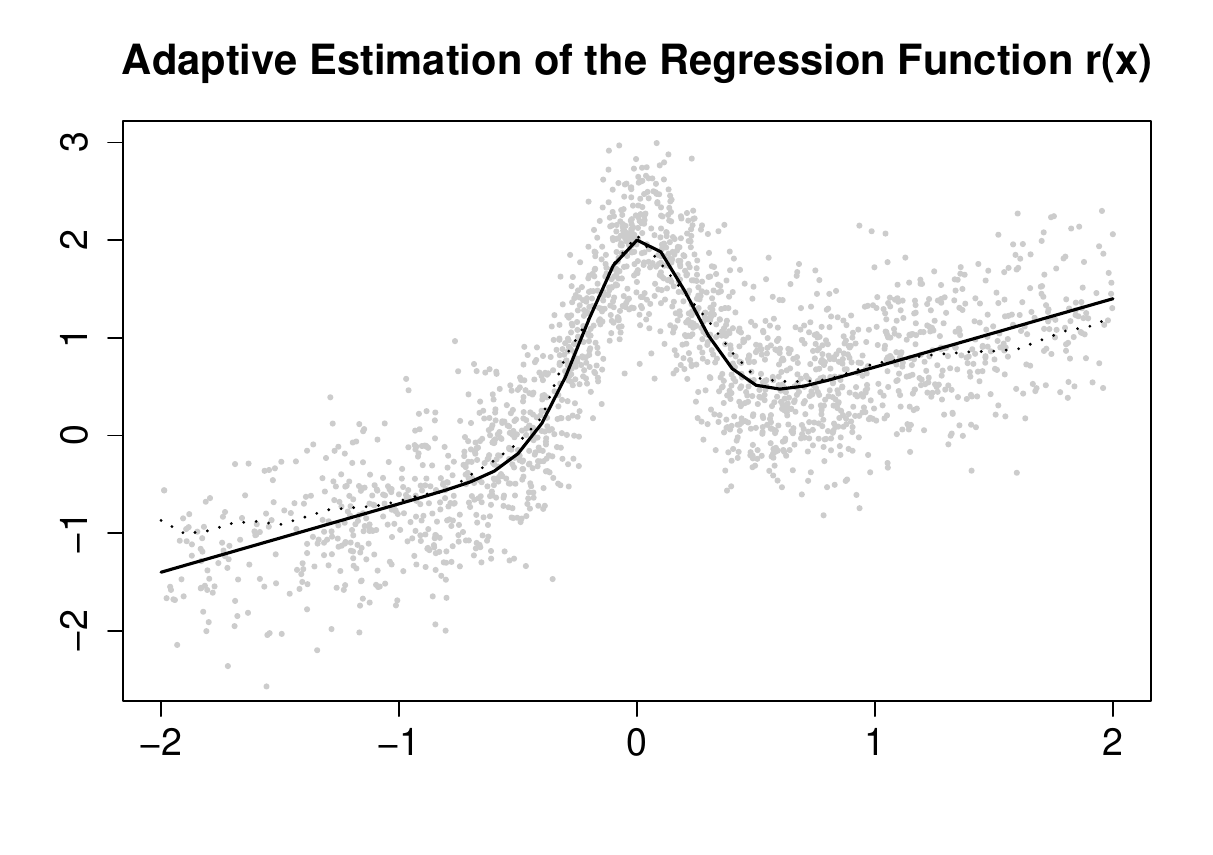}  
	\end{subfigure}
	\begin{subfigure}
	\centering
	\includegraphics[width=7.66cm]{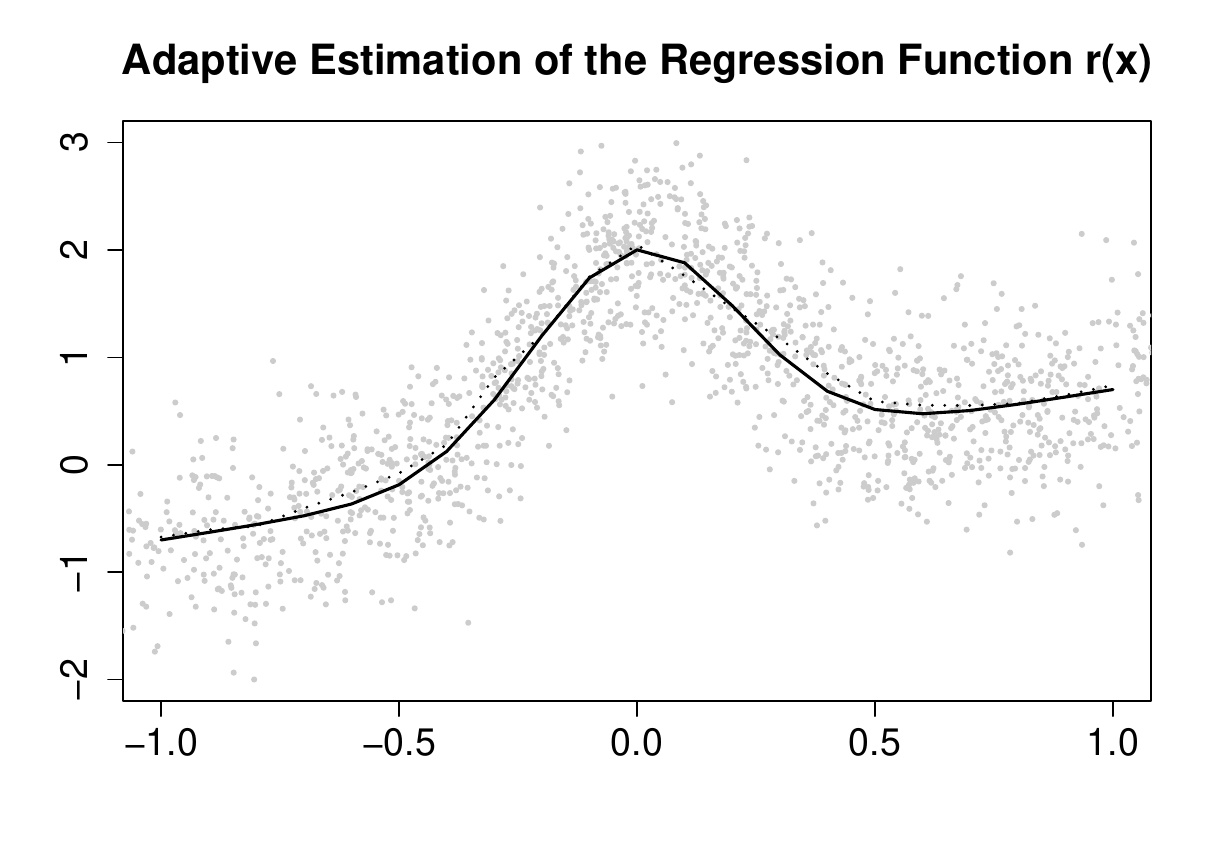}  
	\end{subfigure}
	\caption{The upper figure displays the scatter plot $\{(X_i ,Y_i )\}_{i=1}^n$. The regression function $r$ is represented by a solid line, while the estimator $\{\hat{r}_{\hat{h}_i/\hat{h}^*_i}(x_i)\}_{i=1}^{41}$ on the interval $[-2,2]$ is represented by dots, where $x_{i}=-2+0.1(i-1)$. The lower figure shows the same representations on the grid $x_{i}=-1+0.1(i-1)$, for $i=1,\cdots,21$, over the interval $[-1,1]$. Both estimators were calibrated at $\gamma_{g}=\gamma_6^*$ and $\gamma_{m}=\gamma_2$. In addition, $n=2000$, $q=100$, and $\sigma=0.5$ were used.}
    \label{Crchu-22rchu-11}
\end{figure}

As can be seen from Figure~\ref{CHistyg}, the graph of the density function $g$ does not cover the entire interval $[-2,2]$. This is because the function was evaluated on the random grid $\{\Tilde{x}_i\}_{i=1}^q$, and in this particular case only a few observations were sampled near the endpoints of the interval. This observation highlights potential difficulties in estimating the regression function $r$ near the boundaries of the interval $[-2,2]$.

In Figure~\ref{CCalibraciong}, the left panel shows that $Error_{g}(\gamma^{*})$ attains its minimum at $\gamma_{6}^*$. Moreover, the right panel displays the density function $g$ as a solid line and the estimator $\hat{g}_{\hat{h}^*_i}(x_{i})$, calibrated at $\gamma^*_6$, as points evaluated on the grid $\{x_i\}_{i=1}^{41}$ over the interval $[-2,2]$, where $x_{i}=-2+0.1(i-1)$. In both panels, $n=2000$, $q=100$, and $\sigma=0.5$ were used.

In Figure~\ref{CCalibracionm}, the left panel shows that $Error_{m}(\gamma)$ attains its minimum at $\gamma_{2}$. Moreover, the right panel displays the point cloud
$\left\{ \left(\Tilde{x}_i ,\Tilde{y}_{i}\hat{g}_{\hat{h}_i^{*\gamma_g}}(\Tilde{x}_i)\right) \right\}_{i=1}^q$.
The function $m$ is represented by a solid line, while the estimator $\hat{m}_{\hat{h}_i}(x_{i})$, calibrated at $\gamma_2$, is represented by points and evaluated on the grid $\{x_i\}_{i=1}^{41}$ over the interval $[-2,2]$, where $x_{i}=-2+0.1(i-1)$. In both panels, $n=2000$, $q=100$, and $\sigma=0.5$ were used.

In Figure~\ref{Crchu-22rchu-11}, the left panel displays the scatter plot $\{(X_i ,Y_i )\}_{i=1}^n$. The regression function $r$ is represented by a solid line, while the estimator $\{\hat{r}_{\hat{h}_i/\hat{h}^*_i}(x_i)\}_{i=1}^{41}$ is represented by points on the interval $[-2,2]$, where $x_{i}=-2+0.1(i-1)$. The right panel displays the scatter plot $\{(X_i ,Y_i )\}_{i=1}^n$. The regression function $r$ is represented by a solid line, while the estimator $\{\hat{r}_{\hat{h}_i/\hat{h}^*_i}(x_i)\}_{i=1}^{21}$ is represented by points on the interval $[-1,1]$, where $x_{i}=-1+0.1(i-1)$. In both cases, $n=2000$, $q=100$, and $\sigma=0.5$ were used.

It can be observed that the regression function estimator exhibits poorer performance near the boundaries of the interval $[-2,2]$. For this reason, the final estimation is restricted to the interval $[-1,1]$. Consequently, both the final estimator and the values of $\widehat{MSE}$ and $\widehat{MISE}$ reported in Sections~\ref{CEstGL} and~\ref{CEstErrores} were computed over this interval.

\bibliography{sample}

\begin{thebibliography}{22}
\providecommand{\natexlab}[1]{#1}
\providecommand{\url}[1]{\texttt{#1}}
\expandafter\ifx\csname urlstyle\endcsname\relax
  \providecommand{\doi}[1]{doi: #1}\else
  \providecommand{\doi}{doi: \begingroup \urlstyle{rm}\Url}\fi

\bibitem[Asin and Johannes(2017)]{asin2017adaptive}
N.~Asin and J.~Johannes.
\newblock Adaptive nonparametric estimation in the presence of dependence.
\newblock \emph{Journal of Nonparametric Statistics}, 29\penalty0 (4):\penalty0
  694--730, 2017.

\bibitem[Baraud et~al.(2001)Baraud, Comte, and Viennet]{baraud2001model}
Y.~Baraud, F.~Comte, and G.~Viennet.
\newblock Model selection for (auto-) regression with dependent data.
\newblock \emph{ESAIM: Probability and Statistics}, 5:\penalty0 33--49, 2001.

\bibitem[Bertin and Klutchnikoff(2017)]{bertin2017pointwise}
K.~Bertin and N.~Klutchnikoff.
\newblock Pointwise adaptive estimation of the marginal density of a weakly
  dependent process.
\newblock \emph{Journal of Statistical Planning and Inference}, 187:\penalty0
  115--129, 2017.

\bibitem[Bertin et~al.(2016)Bertin, Lacour, and Rivoirard]{BLR}
K.~Bertin, C.~Lacour, and V.~Rivoirard.
\newblock Adaptive pointwise estimation of conditional density function.
\newblock \emph{Annales de l'Institut Henri Poincar{\'e}, Probabilit{\'e}s et
  Statistiques}, 52\penalty0 (2):\penalty0 939--980, 2016.

\bibitem[Bertin et~al.(2020{\natexlab{a}})Bertin, Klutchnikoff, L{\'e}on, and
  Prieur]{bertin2020adaptivebis}
K.~Bertin, N.~Klutchnikoff, J.~R. L{\'e}on, and C.~Prieur.
\newblock Adaptive density estimation on bounded domains under mixing
  conditions.
\newblock \emph{Electronic Journal of Statistics}, 14\penalty0 (1):\penalty0
  2198--2237, 2020{\natexlab{a}}.

\bibitem[Bertin et~al.(2020{\natexlab{b}})Bertin, Klutchnikoff, Panloup, and
  Varvenne]{bertin2020adaptive}
K.~Bertin, N.~Klutchnikoff, F.~Panloup, and M.~Varvenne.
\newblock Adaptive estimation of the stationary density of a stochastic
  differential equation driven by a fractional brownian motion.
\newblock \emph{Statistical inference for stochastic processes},
  2020{\natexlab{b}}.

\bibitem[Bertin et~al.(2026)Bertin, Fermin, and
  Padrino]{mpadrino2024adaptiregre}
K.~Bertin, L.~Fermin, and M.~Padrino.
\newblock Adaptive estimation in regression models for weakly dependent data
  and explanatory variable with known density.
\newblock \emph{Statistics}, 60\penalty0 (4):\penalty0 1033--1043, 2026.
\newblock \doi{10.1080/02331888.2025.2591680}.

\bibitem[Chichignoud et~al.(2017)Chichignoud, Hoang, Ngoc, and
  Rivoirard]{chichignoud2017adaptive}
M.~Chichignoud, V.~H. Hoang, T.~M.~P. Ngoc, and V.~Rivoirard.
\newblock Adaptive wavelet multivariate regression with errors in variables.
\newblock \emph{Electronic journal of statistics}, 11\penalty0 (1):\penalty0
  682--724, 2017.

\bibitem[Comte and Marie(2021)]{comte}
F.~Comte and N.~Marie.
\newblock {On a Nadaraya-Watson estimator with two bandwidths}.
\newblock \emph{Electronic Journal of Statistics}, 15\penalty0 (1):\penalty0
  2566 -- 2607, 2021.

\bibitem[Comte et~al.(2017)Comte, Prieur, and Samson]{comte2017adaptive}
F.~Comte, C.~Prieur, and A.~Samson.
\newblock Adaptive estimation for stochastic damping hamiltonian systems under
  partial observation.
\newblock \emph{Stochastic Processes and their Applications}, 2017.

\bibitem[Doukhan(1994)]{doukhan1994mixing}
P.~Doukhan.
\newblock Mixing, volume 85 of lecture notes in statistics, 1994.

\bibitem[Doukhan and Neumann(2007)]{doukhan2007probability}
P.~Doukhan and M.~H. Neumann.
\newblock Probability and moment inequalities for sums of weakly dependent
  random variables, with applications.
\newblock \emph{Stochastic Processes and their Applications}, 117\penalty0
  (7):\penalty0 878--903, 2007.

\bibitem[Fermin et~al.(2017)Fermin, Rios, and Rodriguez]{fermin2017robbins}
L.~J. Fermin, R.~Rios, and L.~A. Rodriguez.
\newblock A robbins--monro algorithm for non-parametric estimation of nar
  process with markov switching: Consistency.
\newblock \emph{Journal of Time Series Analysis}, 2017.

\bibitem[Goldenshluger and Lepski(2011)]{goldenshluger2011bandwidth}
A.~Goldenshluger and O.~Lepski.
\newblock Bandwidth selection in kernel density estimation: oracle inequalities
  and adaptive minimax optimality.
\newblock \emph{The Annals of Statistics}, 39\penalty0 (3):\penalty0
  1608--1632, 2011.

\bibitem[Goldenshluger and Lepski(2013)]{GL5}
A.~Goldenshluger and O.~Lepski.
\newblock General selection rule from a family of linear estimators.
\newblock \emph{Theory Prob. Appl.}, 57\penalty0 (2):\penalty0 209--226, 2013.

\bibitem[Goldenshluger and Lepski(2014)]{goldenshluger2014adaptive}
A.~Goldenshluger and O.~Lepski.
\newblock On adaptive minimax density estimation on r\^{}d.
\newblock \emph{Probability Theory and Related Fields}, 159\penalty0
  (3-4):\penalty0 479--543, 2014.

\bibitem[Kouakou and Yod{\'e}(2020)]{Kouakou2020Lp}
K.~F. Kouakou and A.~F.~E. Yod{\'e}.
\newblock $\mathbb{L}_p$-adaptive estimation under partially linear constraint
  in regression model.
\newblock \emph{Journal of Mathematics Research}, 12\penalty0 (6):\penalty0
  74--92, 2020.

\bibitem[Kouakou and Yod{\'e}(2023)]{Kouakou2023Sup}
K.~F. Kouakou and A.~F.~E. Yod{\'e}.
\newblock A sup-norm oracle inequality for a partially linear regression model.
\newblock \emph{Journal of Statistical Planning and Inference}, 222:\penalty0
  132--148, 2023.

\bibitem[Lacour and Massart(2016)]{lacour2016minimalpen}
C.~Lacour and P.~Massart.
\newblock Minimal penalty for goldenshluger-lepski method.
\newblock \emph{Stochastic Processes and their Applications}, 126\penalty0
  (12):\penalty0 3774--3789, 2016.
\newblock \doi{10.1016/j.spa.2016.04.015}.

\bibitem[Lepski and Serdyukova(2014)]{lepski2014adaptive}
O.~Lepski and N.~Serdyukova.
\newblock Adaptive estimation under single-index constraint in a regression
  model.
\newblock \emph{The Annals of Statistics}, 42\penalty0 (1):\penalty0 1--28,
  2014.

\bibitem[Lepski(1990)]{lepski1990problem}
O.~V. Lepski.
\newblock On a problem of adaptive estimation in {G}aussian white noise.
\newblock \emph{Theory of Probability \& Its Applications}, 35\penalty0
  (3):\penalty0 454--466, 1990.

\bibitem[Nadaraya(1989)]{NW}
E.~A. Nadaraya.
\newblock \emph{Nonparametric estimation of probability densities and
  regression curves}.
\newblock Springer, 1989.

\end{thebibliography}

\end{document}